

\ifx\shlhetal\undefinedcontrolsequence\let\shlhetal\relax\fi
\input amstex
\expandafter\ifx\csname mathdefs.tex\endcsname\relax
  \expandafter\gdef\csname mathdefs.tex\endcsname{}
\else \message{Hey!  Apparently you were trying to
  \string\input{mathdefs.tex} twice.   This does not make sense.} 
\errmessage{Please edit your file (probably \jobname.tex) and remove
any duplicate ``\string\input'' lines}\endinput\fi




\catcode`\X=12\catcode`\@=11

\def\n@wcount{\alloc@0\count\countdef\insc@unt}
\def\n@wwrite{\alloc@7\write\chardef\sixt@@n}
\def\n@wread{\alloc@6\read\chardef\sixt@@n}
\def\r@s@t{\relax}\def\v@idline{\par}\def\@mputate#1/{#1}
\def\l@c@l#1X{\firstpart.#1}\def\gl@b@l#1X{#1}\def\t@d@l#1X{{}}

\def\crossrefs#1{\ifx\all#1\let\tr@ce=\all\else\def\tr@ce{#1,}\fi
   \n@wwrite\cit@tionsout\openout\cit@tionsout=\jobname.cit 
   \write\cit@tionsout{\tr@ce}\expandafter\setfl@gs\tr@ce,}
\def\setfl@gs#1,{\def\@{#1}\ifx\@\empty\let\next=\relax
   \else\let\next=\setfl@gs\expandafter\xdef
   \csname#1tr@cetrue\endcsname{}\fi\next}
\def\m@ketag#1#2{\expandafter\n@wcount\csname#2tagno\endcsname
     \csname#2tagno\endcsname=0\let\tail=\all\xdef\all{\tail#2,}
   \ifx#1\l@c@l\let\tail=\r@s@t\xdef\r@s@t{\csname#2tagno\endcsname=0\tail}\fi
   \expandafter\gdef\csname#2cite\endcsname##1{\expandafter
     \ifx\csname#2tag##1\endcsname\relax?\else\csname#2tag##1\endcsname\fi
     \expandafter\ifx\csname#2tr@cetrue\endcsname\relax\else
     \write\cit@tionsout{#2tag ##1 cited on page \folio.}\fi}
   \expandafter\gdef\csname#2page\endcsname##1{\expandafter
     \ifx\csname#2page##1\endcsname\relax?\else\csname#2page##1\endcsname\fi
     \expandafter\ifx\csname#2tr@cetrue\endcsname\relax\else
     \write\cit@tionsout{#2tag ##1 cited on page \folio.}\fi}
   \expandafter\gdef\csname#2tag\endcsname##1{\expandafter
      \ifx\csname#2check##1\endcsname\relax
      \expandafter\xdef\csname#2check##1\endcsname{}%
      \else\immediate\write16{Warning: #2tag ##1 used more than once.}\fi
      \multit@g{#1}{#2}##1/X%
      \write\t@gsout{#2tag ##1 assigned number \csname#2tag##1\endcsname\space
      on page \number\count0.}%
   \csname#2tag##1\endcsname}}

\def\multit@g#1#2#3/#4X{\def\t@mp{#4}\ifx\t@mp\empty%
      \global\advance\csname#2tagno\endcsname by 1 
      \expandafter\xdef\csname#2tag#3\endcsname
      {#1\number\csname#2tagno\endcsnameX}%
   \else\expandafter\ifx\csname#2last#3\endcsname\relax
      \expandafter\n@wcount\csname#2last#3\endcsname
      \global\advance\csname#2tagno\endcsname by 1 
      \expandafter\xdef\csname#2tag#3\endcsname
      {#1\number\csname#2tagno\endcsnameX}
      \write\t@gsout{#2tag #3 assigned number \csname#2tag#3\endcsname\space
      on page \number\count0.}\fi
   \global\advance\csname#2last#3\endcsname by 1
   \def\t@mp{\expandafter\xdef\csname#2tag#3/}%
   \expandafter\t@mp\@mputate#4\endcsname
   {\csname#2tag#3\endcsname\lastpart{\csname#2last#3\endcsname}}\fi}
\def\t@gs#1{\def\all{}\m@ketag#1e\m@ketag#1s\m@ketag\t@d@l p
\let\realscite\scite
\let\realstag\stag
   \m@ketag\gl@b@l r \n@wread\t@gsin
   \openin\t@gsin=\jobname.tgs \re@der \closein\t@gsin
   \n@wwrite\t@gsout\openout\t@gsout=\jobname.tgs }
\outer\def\localtags{\t@gs\l@c@l}
\outer\def\globaltags{\t@gs\gl@b@l}
\outer\def\newlocaltag#1{\m@ketag\l@c@l{#1}}
\outer\def\newglobaltag#1{\m@ketag\gl@b@l{#1}}

\newif\ifpr@ 
\def\m@kecs #1tag #2 assigned number #3 on page #4.%
   {\expandafter\gdef\csname#1tag#2\endcsname{#3}
   \expandafter\gdef\csname#1page#2\endcsname{#4}
   \ifpr@\expandafter\xdef\csname#1check#2\endcsname{}\fi}
\def\re@der{\ifeof\t@gsin\let\next=\relax\else
   \read\t@gsin to\t@gline\ifx\t@gline\v@idline\else
   \expandafter\m@kecs \t@gline\fi\let \next=\re@der\fi\next}
\def\pretags#1{\pr@true\pret@gs#1,,}
\def\pret@gs#1,{\def\@{#1}\ifx\@\empty\let\n@xtfile=\relax
   \else\let\n@xtfile=\pret@gs \openin\t@gsin=#1.tgs \message{#1} \re@der 
   \closein\t@gsin\fi \n@xtfile}

\newcount\sectno\sectno=0\newcount\subsectno\subsectno=0
\newif\ifultr@local \def\ultralocal{\ultr@localtrue}
\def\firstpart{\number\sectno}
\def\lastpart#1{\ifcase#1 \or a\or b\or c\or d\or e\or f\or g\or h\or 
   i\or k\or l\or m\or n\or o\or p\or q\or r\or s\or t\or u\or v\or w\or 
   x\or y\or z \fi}

\def\resetall{\global\advance\sectno by 1\subsectno=0
   \gdef\firstpart{\number\sectno}\r@s@t}
\def\resetsub{\global\advance\subsectno by 1
   \gdef\firstpart{\number\sectno.\number\subsectno}\r@s@t}
\def\newsection#1\par{\resetall\vskip0pt plus.3\vsize\penalty-250
   \vskip0pt plus-.3\vsize\bigskip\bigskip
   \message{#1}\leftline{\bf#1}\nobreak\bigskip}
\def\subsection#1\par{\ifultr@local\resetsub\fi
   \vskip0pt plus.2\vsize\penalty-250\vskip0pt plus-.2\vsize
   \bigskip\smallskip\message{#1}\leftline{\bf#1}\nobreak\medskip}


\newdimen\marginshift

\newdimen\margindelta
\newdimen\marginmax
\newdimen\marginmin

\def\margininit{       
\marginmax=3 true cm                  
				      
\margindelta=0.1 true cm              
\marginmin=0.1true cm                 
\marginshift=\marginmin
}    

\def\t@gsjj#1,{\def\@{#1}\ifx\@\empty\let\next=\relax\else\let\next=\t@gsjj
   \def\@@{p}\ifx\@\@@\else
   \expandafter\gdef\csname#1cite\endcsname##1{\citejj{##1}}
   \expandafter\gdef\csname#1page\endcsname##1{?}
   \expandafter\gdef\csname#1tag\endcsname##1{\tagjj{##1}}\fi\fi\next}
\newif\ifshowstuffinmargin
\showstuffinmarginfalse
\def\jjtags{\showstuffinmargintrue
\ifx\all\relax\else\expandafter\t@gsjj\all,\fi}

\def\tagjj#1{\realstag{#1}\mginpar{\zeigen{#1}}}
\def\citejj#1{\zeigen{#1}\mginpar{\rechnen{#1}}}

\def\rechnen#1{\expandafter\ifx\csname stag#1\endcsname\relax ??\else
                           \csname stag#1\endcsname\fi}

\newdimen\theight

\def\marginfont{\sevenrm}

\def\trymarginbox#1{\setbox0=\hbox{\marginfont\hskip\marginshift #1}%
		\global\marginshift\wd0 
		\global\advance\marginshift\margindelta}

\def \mginpar#1{%
\ifvmode\setbox0\hbox to \hsize{\hfill\rlap{\marginfont\quad#1}}%
\ht0 0cm
\dp0 0cm
\box0\vskip-\baselineskip
\else 
             \vadjust{\trymarginbox{#1}%
		\ifdim\marginshift>\marginmax \global\marginshift\marginmin
			\trymarginbox{#1}%
                \fi
             \theight=\ht0
             \advance\theight by \dp0    \advance\theight by \lineskip
             \kern -\theight \vbox to \theight{\rightline{\rlap{\box0}}%
\vss}}\fi}


\def\t@gsoff#1,{\def\@{#1}\ifx\@\empty\let\next=\relax\else\let\next=\t@gsoff
   \def\@@{p}\ifx\@\@@\else
   \expandafter\gdef\csname#1cite\endcsname##1{\zeigen{##1}}
   \expandafter\gdef\csname#1page\endcsname##1{?}
   \expandafter\gdef\csname#1tag\endcsname##1{\zeigen{##1}}\fi\fi\next}
\def\verbatimtags{\showstuffinmarginfalse
\ifx\all\relax\else\expandafter\t@gsoff\all,\fi}
\def\zeigen#1{\hbox{$\langle$}#1\hbox{$\rangle$}}

\def\(#1){\edef\dot@g{\ifmmode\ifinner(\hbox{\noexpand\etag{#1}})
   \else\noexpand\eqno(\hbox{\noexpand\etag{#1}})\fi
   \else(\noexpand\ecite{#1})\fi}\dot@g}

\newif\ifbr@ck
\def\eat#1{}
\def\[#1]{\br@cktrue[\br@cket#1'X]}
\def\br@cket#1'#2X{\def\temp{#2}\ifx\temp\empty\let\next\eat
   \else\let\next\br@cket\fi
   \ifbr@ck\br@ckfalse\br@ck@t#1,X\else\br@cktrue#1\fi\next#2X}
\def\br@ck@t#1,#2X{\def\temp{#2}\ifx\temp\empty\let\neext\eat
   \else\let\neext\br@ck@t\def\temp{,}\fi
   \def\teemp{#1}\ifx\teemp\empty\else\rcite{#1}\fi\temp\neext#2X}
\def\resetbr@cket{\gdef\[##1]{[\rtag{##1}]}}
\def\references{\resetbr@cket\newsection References\par}

\newtoks\symb@ls\newtoks\s@mb@ls\newtoks\p@gelist\n@wcount\ftn@mber
    \ftn@mber=1\newif\ifftn@mbers\ftn@mbersfalse\newif\ifbyp@ge\byp@gefalse
\def\defm@rk{\ifftn@mbers\n@mberm@rk\else\symb@lm@rk\fi}
\def\n@mberm@rk{\xdef\m@rk{{\the\ftn@mber}}%
    \global\advance\ftn@mber by 1 }
\def\rot@te#1{\let\temp=#1\global#1=\expandafter\r@t@te\the\temp,X}
\def\r@t@te#1,#2X{{#2#1}\xdef\m@rk{{#1}}}
\def\b@@st#1{{$^{#1}$}}\def\str@p#1{#1}
\def\symb@lm@rk{\ifbyp@ge\rot@te\p@gelist\ifnum\expandafter\str@p\m@rk=1 
    \s@mb@ls=\symb@ls\fi\write\f@nsout{\number\count0}\fi \rot@te\s@mb@ls}
\def\byp@ge{\byp@getrue\n@wwrite\f@nsin\openin\f@nsin=\jobname.fns 
    \n@wcount\currentp@ge\currentp@ge=0\p@gelist={0}
    \re@dfns\closein\f@nsin\rot@te\p@gelist
    \n@wread\f@nsout\openout\f@nsout=\jobname.fns }
\def\m@kelist#1X#2{{#1,#2}}
\def\re@dfns{\ifeof\f@nsin\let\next=\relax\else\read\f@nsin to \f@nline
    \ifx\f@nline\v@idline\else\let\t@mplist=\p@gelist
    \ifnum\currentp@ge=\f@nline
    \global\p@gelist=\expandafter\m@kelist\the\t@mplistX0
    \else\currentp@ge=\f@nline
    \global\p@gelist=\expandafter\m@kelist\the\t@mplistX1\fi\fi
    \let\next=\re@dfns\fi\next}
\def\symbols#1{\symb@ls={#1}\s@mb@ls=\symb@ls} 
\def\bigsymbol{\textstyle}
\symbols{\bigsymbol\ast,\dagger,\ddagger,\sharp,\flat,\natural,\star}
\def\ftnumbers{\ftn@mberstrue} \def\ftsymbols{\ftn@mbersfalse}
\def\paginal{\byp@ge} \def\resetftnumbers{\ftn@mber=1}
\def\ftnote#1{\defm@rk\expandafter\expandafter\expandafter\footnote
    \expandafter\b@@st\m@rk{#1}}

\long\def\jump#1\endjump{}
\def\ssum{\mathop{\lower .1em\hbox{$\textstyle\Sigma$}}\nolimits}

\def\qed{\nobreak\kern 1em \vrule height .5em width .5em depth 0em}
\def\newneq{\hbox{\rlap{\hbox to 1\wd9{\hss$=$\hss}}\raise .1em 
   \hbox to 1\wd9{\hss$\scriptscriptstyle/$\hss}}}
\def\subsetne{\setbox9 = \hbox{$\subset$}\mathrel{\hbox{\rlap
   {\lower .4em \newneq}\raise .13em \hbox{$\subset$}}}}
\def\supsetne{\setbox9 = \hbox{$\subset$}\mathrel{\hbox{\rlap
   {\lower .4em \newneq}\raise .13em \hbox{$\supset$}}}}

\def\vbar{\mathchoice{\vrule height6.3ptdepth-.5ptwidth.8pt\kern-.8pt}
   {\vrule height6.3ptdepth-.5ptwidth.8pt\kern-.8pt}
   {\vrule height4.1ptdepth-.35ptwidth.6pt\kern-.6pt}
   {\vrule height3.1ptdepth-.25ptwidth.5pt\kern-.5pt}}
\def\f@dge{\mathchoice{}{}{\mkern.5mu}{\mkern.8mu}}
\def\b@c#1#2{{\rm \mkern#2mu\vbar\mkern-#2mu#1}}
\def\b@b#1{{\rm I\mkern-3.5mu #1}}
\def\b@a#1#2{{\rm #1\mkern-#2mu\f@dge #1}}
\def\bb#1{{\count4=`#1 \advance\count4by-64 \ifcase\count4\or\b@a A{11.5}\or
   \b@b B\or\b@c C{5}\or\b@b D\or\b@b E\or\b@b F \or\b@c G{5}\or\b@b H\or
   \b@b I\or\b@c J{3}\or\b@b K\or\b@b L \or\b@b M\or\b@b N\or\b@c O{5} \or
   \b@b P\or\b@c Q{5}\or\b@b R\or\b@a S{8}\or\b@a T{10.5}\or\b@c U{5}\or
   \b@a V{12}\or\b@a W{16.5}\or\b@a X{11}\or\b@a Y{11.7}\or\b@a Z{7.5}\fi}}

\catcode`\X=11 \catcode`\@=12


\expandafter\ifx\csname citeadd.tex\endcsname\relax
\expandafter\gdef\csname citeadd.tex\endcsname{}
\else \message{Hey!  Apparently you were trying to
\string\input{citeadd.tex} twice.   This does not make sense.} 
\errmessage{Please edit your file (probably \jobname.tex) and remove
any duplicate ``\string\input'' lines}\endinput\fi

\sectno=-1   
\localtags
\NoBlackBoxes
\define\mr{\medskip\roster}
\define\sn{\smallskip\noindent}
\define\mn{\medskip\noindent}
\define\bn{\bigskip\noindent}
\define\ub{\underbar}
\define\wilog{\text{without loss of generality}}
\define\ermn{\endroster\medskip\noindent}
\define\dbca{\dsize\bigcap}
\define\dbcu{\dsize\bigcup}
\define \nl{\newline}
\magnification=\magstep 1
\documentstyle{amsppt}

{    
\catcode`@11

\ifx\alicetwothousandloaded@\relax
  \endinput\else\global\let\alicetwothousandloaded@\relax\fi

\gdef\subjclass{\let\savedef@\subjclass
 \def\subjclass##1\endsubjclass{\let\subjclass\savedef@
   \toks@{\def\usualspace{{\rm\enspace}}\eightpoint}%
   \toks@@{##1\unskip.}%
   \edef\thesubjclass@{\the\toks@
     \frills@{{\noexpand\rm2000 {\noexpand\it Mathematics Subject
       Classification}.\noexpand\enspace}}%
     \the\toks@@}}%
  \nofrillscheck\subjclass}
} 

\pageheight{8.5truein}
\topmatter
\title{On nice equivalence relations on ${}^\lambda 2$} \endtitle
\author {Saharon Shelah \thanks {\null\newline I would like to thank 
Alice Leonhardt for the beautiful typing \null\newline
This research was supported by The Israel Science Foundation
founded by the Israel Academy of Sciences and Humanities. Publication
724.}\endthanks} \endauthor    

\affil{Institute of Mathematics\\
 The Hebrew University\\
 Jerusalem, Israel
 \medskip
 Rutgers University\\
 Mathematics Department\\
 New Brunswick, NJ  USA} \endaffil
\endtopmatter
\document  

\expandafter\ifx\csname alice2jlem.tex\endcsname\relax
  \expandafter\xdef\csname alice2jlem.tex\endcsname{\the\catcode`@}
\else \message{Hey!  Apparently you were trying to
\string\input{alice2jlem.tex}  twice.   This does not make sense.}
\errmessage{Please edit your file (probably \jobname.tex) and remove
any duplicate ``\string\input'' lines}\endinput\fi

\expandafter\ifx\csname bib4plain.tex\endcsname\relax
  \expandafter\gdef\csname bib4plain.tex\endcsname{}
\else \message{Hey!  Apparently you were trying to \string\input
  bib4plain.tex twice.   This does not make sense.}
\errmessage{Please edit your file (probably \jobname.tex) and remove
any duplicate ``\string\input'' lines}\endinput\fi

\def\renewcommand{\newcommand}	       
\edef\cite{\the\catcode`@}%
\catcode`@ = 11
\let\@oldatcatcode = \cite
\chardef\@letter = 11
\chardef\@other = 12
%
%
%
%
\def\@innerdef#1#2{\edef#1{\expandafter\noexpand\csname #2\endcsname}}%
%
%
\@innerdef\@innernewcount{newcount}%
\@innerdef\@innernewdimen{newdimen}%
\@innerdef\@innernewif{newif}%
\@innerdef\@innernewwrite{newwrite}%
%
%
%
\def\@gobble#1{}%
%
%
%
\ifx\inputlineno\@undefined
   \let\@linenumber = \empty 
\else
   \def\@linenumber{\the\inputlineno:\space}%
\fi
%
%
%
\def\@futurenonspacelet#1{\def\cs{#1}%
   \afterassignment\@stepone\let\@nexttoken=
}%
\begingroup 
\def\\{\global\let\@stoken= }%
\\ 
\endgroup
\def\@stepone{\expandafter\futurelet\cs\@steptwo}%
\def\@steptwo{\expandafter\ifx\cs\@stoken\let\@@next=\@stepthree
   \else\let\@@next=\@nexttoken\fi \@@next}%
\def\@stepthree{\afterassignment\@stepone\let\@@next= }%
%
%
%
\def\@getoptionalarg#1{%
   \let\@optionaltemp = #1%
   \let\@optionalnext = \relax
   \@futurenonspacelet\@optionalnext\@bracketcheck
}%
%
%
\def\@bracketcheck{%
   \ifx [\@optionalnext
      \expandafter\@@getoptionalarg
   \else
      \let\@optionalarg = \empty
      \expandafter\@optionaltemp
   \fi
}%
\def\@@getoptionalarg[#1]{%
   \def\@optionalarg{#1}%
   \@optionaltemp
}%
%
%
%
\def\@nnil{\@nil}%
\def\@fornoop#1\@@#2#3{}%
\def\@for#1:=#2\do#3{%
   \edef\@fortmp{#2}%
   \ifx\@fortmp\empty \else
      \expandafter\@forloop#2,\@nil,\@nil\@@#1{#3}%
   \fi
}%
\def\@forloop#1,#2,#3\@@#4#5{\def#4{#1}\ifx #4\@nnil \else
       #5\def#4{#2}\ifx #4\@nnil \else#5\@iforloop #3\@@#4{#5}\fi\fi
}%
\def\@iforloop#1,#2\@@#3#4{\def#3{#1}\ifx #3\@nnil
       \let\@nextwhile=\@fornoop \else
      #4\relax\let\@nextwhile=\@iforloop\fi\@nextwhile#2\@@#3{#4}%
}%
%
%
%
\@innernewif\if@fileexists
\def\@testfileexistence{\@getoptionalarg\@finishtestfileexistence}%
\def\@finishtestfileexistence#1{%
   \begingroup
      \def\extension{#1}%
      \immediate\openin0 =
         \ifx\@optionalarg\empty\jobname\else\@optionalarg\fi
         \ifx\extension\empty \else .#1\fi
         \space
      \ifeof 0
         \global\@fileexistsfalse
      \else
         \global\@fileexiststrue
      \fi
      \immediate\closein0
   \endgroup
}%
%
%
%
%
\def\bibliographystyle#1{%
   \@readauxfile
   \@writeaux{\string\bibstyle{#1}}%
}%
\let\bibstyle = \@gobble
%
%
\let\bblfilebasename = \jobname
\def\bibliography#1{%
   \@readauxfile
   \@writeaux{\string\bibdata{#1}}%
   \@testfileexistence[\bblfilebasename]{bbl}%
   \if@fileexists
      \nobreak
      \@readbblfile
   \fi
}%
\let\bibdata = \@gobble
%
%
\def\nocite#1{%
   \@readauxfile
   \@writeaux{\string\citation{#1}}%
}%
\@innernewif\if@notfirstcitation
%
%
\def\cite{\@getoptionalarg\@cite}%
%
%
\def\@cite#1{%
   \let\@citenotetext = \@optionalarg
   \printcitestart
   \nocite{#1}%
   \@notfirstcitationfalse
   \@for \@citation :=#1\do
   {%
      \expandafter\@onecitation\@citation\@@
   }%
   \ifx\empty\@citenotetext\else
      \printcitenote{\@citenotetext}%
   \fi
   \printcitefinish
}%
\def\@onecitation#1\@@{%
   \if@notfirstcitation
      \printbetweencitations
   \fi
   \expandafter \ifx \csname\@citelabel{#1}\endcsname \relax
      \if@citewarning
         \message{\@linenumber Undefined citation `#1'.}%
      \fi
      \expandafter\gdef\csname\@citelabel{#1}\endcsname{%
\strut
\vadjust{\vskip-\dp\strutbox
\vbox to 0pt{\vss\parindent0cm \leftskip=\hsize 
\advance\leftskip3mm
\advance\hsize 4cm\strut\openup-4pt 
\rightskip 0cm plus 1cm minus 0.5cm ?  #1 ?\strut}}
         {\tt
            \escapechar = -1
            \nobreak\hskip0pt
            \expandafter\string\csname#1\endcsname
            \nobreak\hskip0pt
         }%
      }%
   \fi
   \csname\@citelabel{#1}\endcsname
   \@notfirstcitationtrue
}%
%
%
\def\@citelabel#1{b@#1}%
%
%
\def\@citedef#1#2{\expandafter\gdef\csname\@citelabel{#1}\endcsname{#2}}%
%
%
%
\def\@readbblfile{%
   \ifx\@itemnum\@undefined
      \@innernewcount\@itemnum
   \fi
   \begingroup
      \def\begin##1##2{%
         \setbox0 = \hbox{\biblabelcontents{##2}}%
         \biblabelwidth = \wd0
      }%
      \def\end##1{}
      %
      %
      \@itemnum = 0
      \def\bibitem{\@getoptionalarg\@bibitem}%
      \def\@bibitem{%
         \ifx\@optionalarg\empty
            \expandafter\@numberedbibitem
         \else
            \expandafter\@alphabibitem
         \fi
      }%
      \def\@alphabibitem##1{%
         \expandafter \xdef\csname\@citelabel{##1}\endcsname {\@optionalarg}%
         \ifx\biblabelprecontents\@undefined
            \let\biblabelprecontents = \relax
         \fi
         \ifx\biblabelpostcontents\@undefined
            \let\biblabelpostcontents = \hss
         \fi
         \@finishbibitem{##1}%
      }%
      \def\@numberedbibitem##1{%
         \advance\@itemnum by 1
         \expandafter \xdef\csname\@citelabel{##1}\endcsname{\number\@itemnum}%
         \ifx\biblabelprecontents\@undefined
            \let\biblabelprecontents = \hss
         \fi
         \ifx\biblabelpostcontents\@undefined
            \let\biblabelpostcontents = \relax
         \fi
         \@finishbibitem{##1}%
      }%
      \def\@finishbibitem##1{%
         \biblabelprint{\csname\@citelabel{##1}\endcsname}%
         \@writeaux{\string\@citedef{##1}{\csname\@citelabel{##1}\endcsname}}%
         \ignorespaces
      }%
      %
      %
      \let\em = \bblem
      \let\newblock = \bblnewblock
      \let\sc = \bblsc
      \frenchspacing
      \clubpenalty = 4000 \widowpenalty = 4000
      \tolerance = 10000 \hfuzz = .5pt
      \everypar = {\hangindent = \biblabelwidth
                      \advance\hangindent by \biblabelextraspace}%
      \bblrm
      \parskip = 1.5ex plus .5ex minus .5ex
      \biblabelextraspace = .5em
      \bblhook
      \input \bblfilebasename.bbl
   \endgroup
}%
%
%
\@innernewdimen\biblabelwidth
\@innernewdimen\biblabelextraspace
%
%
%
\def\biblabelprint#1{%
   \noindent
   \hbox to \biblabelwidth{%
      \biblabelprecontents
      \biblabelcontents{#1}%
      \biblabelpostcontents
   }%
   \kern\biblabelextraspace
}%
%
%
%
\def\biblabelcontents#1{{\bblrm [#1]}}%
%
%
\def\bblrm{\rm}%
%
%
\def\bblem{\it}%
%
%
\def\bblsc{\ifx\@scfont\@undefined
              \font\@scfont = cmcsc10
           \fi
           \@scfont
}%
%
%
\def\bblnewblock{\hskip .11em plus .33em minus .07em }%
%
%
\let\bblhook = \empty
%
%
%
\def\printcitestart{[}
\def\printcitefinish{]}
\def\printbetweencitations{, }
\def\printcitenote#1{, #1}
%
%
%
\let\citation = \@gobble
%
%
%
\@innernewcount\@numparams
%
%
\def\newcommand#1{%
   \def\@commandname{#1}%
   \@getoptionalarg\@continuenewcommand
}%
%
%
\def\@continuenewcommand{%
   \@numparams = \ifx\@optionalarg\empty 0\else\@optionalarg \fi \relax
   \@newcommand
}%
%
%
\def\@newcommand#1{%
   \def\@startdef{\expandafter\edef\@commandname}%
   \ifnum\@numparams=0
      \let\@paramdef = \empty
   \else
      \ifnum\@numparams>9
         \errmessage{\the\@numparams\space is too many parameters}%
      \else
         \ifnum\@numparams<0
            \errmessage{\the\@numparams\space is too few parameters}%
         \else
            \edef\@paramdef{%
               \ifcase\@numparams
                  \empty  No arguments.
               \or ####1%
               \or ####1####2%
               \or ####1####2####3%
               \or ####1####2####3####4%
               \or ####1####2####3####4####5%
               \or ####1####2####3####4####5####6%
               \or ####1####2####3####4####5####6####7%
               \or ####1####2####3####4####5####6####7####8%
               \or ####1####2####3####4####5####6####7####8####9%
               \fi
            }%
         \fi
      \fi
   \fi
   \expandafter\@startdef\@paramdef{#1}%
}%
%
%
%
%
\def\@readauxfile{%
   \if@auxfiledone \else 
      \global\@auxfiledonetrue
      \@testfileexistence{aux}%
      \if@fileexists
         \begingroup
            \endlinechar = -1
            \catcode`@ = 11
            \input \jobname.aux
         \endgroup
      \else
         \message{\@undefinedmessage}%
         \global\@citewarningfalse
      \fi
      \immediate\openout\@auxfile = \jobname.aux
   \fi
}%
%
%
\newif\if@auxfiledone
\ifx\noauxfile\@undefined \else \@auxfiledonetrue\fi
%
%
%
%
\@innernewwrite\@auxfile
\def\@writeaux#1{\ifx\noauxfile\@undefined \write\@auxfile{#1}\fi}%
%
%
%
\ifx\@undefinedmessage\@undefined
   \def\@undefinedmessage{No .aux file; I won't give you warnings about
                          undefined citations.}%
\fi
%
%
\@innernewif\if@citewarning
\ifx\noauxfile\@undefined \@citewarningtrue\fi
%
%
%
\catcode`@ = \@oldatcatcode


\def\widestnumber#1#2{}

\def\rm{\fam0 \tenrm}

\def\fakesubhead#1\endsubhead{\bigskip\noindent{\bf#1}\par}



%
%
%

%

\font\textrsfs=rsfs10
\font\scriptrsfs=rsfs7
\font\scriptscriptrsfs=rsfs5

\newfam\rsfsfam
\textfont\rsfsfam=\textrsfs
\scriptfont\rsfsfam=\scriptrsfs
\scriptscriptfont\rsfsfam=\scriptscriptrsfs

\edef\oldcatcodeofat{\the\catcode`\@}
\catcode`\@11

\def\Cal@@#1{\noaccents@ \fam \rsfsfam #1}

\catcode`\@\oldcatcodeofat


\expandafter\ifx \csname margininit\endcsname \relax\else\margininit\fi

\bn

\head {\S0 Introduction} \endhead  \resetall \sectno=0
\bigskip

The main question here is the possible generalization of the following
theorem on ``simple" equivalence relation on ${}^\omega 2$ to higher
cardinals.
\proclaim{\stag{0.1} Theorem}  1) Assume that
\mr
\item "{$(a)$}"  $E$ is a Borel 2-place relation on ${}^\omega 2$
\sn
\item "{$(b)$}"  $E$ is an equivalence relation
\sn
\item "{$(c)$}"  if $\eta,\nu \in {}^\omega 2$ and
$(\exists!n)(\eta(n) \ne \nu(n))$, \ub{then} $\eta,\nu$ are not
$E$-equivalent.
\ermn
\ub{Then} there is a perfect subset of ${}^\omega 2$ of pairwise non
$E$-equivalent members. \nl
2) Instead of ``$E$ is Borel", ``$E$ is analytic (or even a Borel
combination of analytic relations)" is enough. \nl
3) If $E$ is a $\pi^1_2$ relation which is an equivalence relation
satisfying clauses (b) + (c) in $\bold V^{\text{Cohen}}$, \ub{then} the
conclusion of (1) holds.
\endproclaim
\bn
See \cite{Sh:273}, it was used to prove a result on the homotopy
group: if $X$ is Hausdorff metric, compact, separable, arc-connected,
locally arc-connected and the homotopy group is not finitely generated
then it has cardinality continuum; the proof of \scite{0.1} used
forcing in \cite{Sh:273}, see \cite{PaSr98} for eliminating the
forcing. 

We may restrict $E$ to be like $r_p(\text{Ext}(G,\Bbb Z))$ or just
closer to group theory as in Grossberg Shelah \cite{GrSh:302}, 
\cite{GrSh:302a}, Mekler Roslanowski Shelah \cite{MRSh:314},
\cite{Sh:664}.  In \S5 we say somewhat more. \nl
We here continue \cite{Sh:664} but do not rely on it.

Turning to ${}^\lambda 2$ the problem split according to the character
of $\lambda$ and the ``simplicity" of $E$.  If $E$ is $\pi^1_1$ and
$\lambda = \lambda^{< \lambda}$ and $\lambda \ge \beth_\omega$ (or
just $(Dl)_\lambda$ holds) a generalization holds.  If $E$ is $\Sigma^1_1,\lambda =
\lambda^{< \lambda}$ the generalization in general fails; all this in
\S1.  Now if $\lambda$ is singular, strong limit for
simplicity, it is natural to consider ${}^{\text{cf}(\lambda)}\lambda$
instead of ${}^\lambda 2$.  If $\lambda$ has uncountable cofinality we
get strong negative results in \S2.  If $\lambda$ has countable cofinality,
and is the limit of ``somewhat large cardinals", e.g. measurable
cardinals, (but $\lambda = \aleph_\omega$ may be O.K., i.e. consistently) the
generalization holds (in \S3), but if the universe is closed to $L$ (e.g. in
$L$ there is no weakly compact or just no Erdos cardinal) then we get
negative results (see \S4).  Note that theorems of the form ``if $E$ has many
equivalence classes it has continuum" do not generalize well, see
\cite{ShVs:719} even for weakly compact.  
\bigskip

\definition{\stag{0.2} Definition}  For a cardinal $\lambda$ and let
${\Cal B}_\lambda$ be ${}^\lambda 2$ (or ${}^\lambda \lambda$ or
${}^{\text{cf}(\lambda)} \lambda)$. \nl
1)  For a logic ${\Cal L}$ we say that $E$ is a ${\Cal L}$-\ub{nice},
say 2-place for simplicity, relation on ${\Cal B}$ \ub{if} there is
a model $M$ with universe $\lambda$ and finite vocabulary $\tau$, and
unary function symbols $F_1,F_2 \notin \tau$ (denoting possibly partial
unary functions), such that letting $\tau^+ = \tau \cup \{F_1,F_2\}$, for some
sentence $\psi = \psi(F_1,F_2)$ in $L(\tau)$ we have
\mr
\item "{$\circledcirc$}"  for any $\eta_1,\eta_2 \in {\Cal B}$ letting
$M_{\eta_1,\eta_2} = (M,\eta_1,\eta_2)$ be the $\tau^+$-model
expanding $M$ with $F^{M_{\eta_1,\eta_2}}_\ell = \eta_\ell$ we have
\nl
$\eta_1 E \eta_2 \Leftrightarrow (M,\eta_1,\eta_2) \models \psi$. \nl
We may write $M \models \psi[\eta_1,\eta_2]$ and
$\psi[\eta_1,\eta_2,M]$ or $\psi(x,y,M)$ or write $a \subseteq
\lambda$ coding $M$ instead of $M$.
\ermn
2) $E$ is a $\pi^1_1$-relation on ${\Cal B}$ means that above we allow
$\psi$ to be of the form $(\forall X)\varphi,\varphi$ first order or
even inductive logic (i.e. we have variables on sets and are allowed
to form the first fix point); $X$
vary on sets.  Similarly $\Sigma^1_1,\pi^1_2$, projective; writing
nice means ${\Cal L}$ is first order + definition by induction.  We
may write nice$_{\Cal B},\Sigma^1_1[{\Cal B}]$ etc, and may replace
${\Cal B}$ by $\lambda$ if ${\Cal B} = {\Cal B}_\lambda$.  We write
very nice for ${\Cal L}$-nice, ${\Cal L}$ first order.
\enddefinition
\bn
\ub{Notation}:

$(\forall^* i < \delta)$ means ``for every large enough $i < \delta$".

$J^{\text{bd}}_\delta$ is the ideal of bounded subsets of $\delta$.
\bigskip

\definition{\stag{0.3} Definition}  Let $(D \ell)_\lambda$ means that $\lambda$
is regular, uncountable and there is $\bar{\Cal P} = \langle {\Cal
P}_\alpha:\alpha < \lambda \rangle$ such that ${\Cal P}_\alpha$ is a
family of $< \lambda$ subsets of $\alpha$ and for every $X \subseteq
\lambda$ the set $\{\delta < \lambda:X \cap \delta \in {\Cal
P}_\delta\}$ is stationary; hence $\lambda = \lambda^{< \lambda}$.  
[By \cite{Sh:460}, $\lambda = \lambda^{< \lambda} 
\ge \beth_\omega \Rightarrow (D \ell)_\lambda$ and (by Kunen)
$\lambda = \mu^+ \Rightarrow (D \ell)_\lambda \equiv \diamondsuit_\lambda$].
\enddefinition
\bigskip

\definition{\stag{0.4} Definition}  ${\Cal Q} \subseteq {}^\lambda 2$
is called perfect or $\lambda$-perfect if:
\mr
\item "{$(a)$}"  $A \ne \emptyset$
\sn
\item "{$(b)$}"  if $\eta \in A$ then $\{\ell g(\eta \cap \nu):\nu \in
A \backslash \{\eta\}\} \subseteq \lambda$ is unbounded
\sn
\item "{$(c)$}"  $\{\eta \restriction \zeta:\eta \in A,\zeta \le
\lambda\}$ is closed under the union of $\triangleleft$-increasing sequences.
\ermn
Equivalently, ${\Cal Q} = \langle \rho_\eta:\eta \in {}^\lambda 2
\rangle$ such that
\mr
\item "{$(a)'$}"  $\rho_\eta \in {}^\lambda 2$
\sn
\item "{$(b)'$}"  $\eta_1 \ne \eta_2 \in {}^\lambda 2 \Rightarrow
\rho_{\eta_2} \ne \rho_{\eta_2}$
\sn
\item "{$(c)'$}"  if $\eta_0,\eta_1,\eta_2 \in {}^\lambda 2$ are
distinct and if $\eta_1 \cap \eta_2 \triangleleft \eta_1 \cap \eta_0$
(so $\eta_1 \cap \eta_2 \ne \eta_1 \cap \eta_0)$ then $\rho_{\eta_1}
\cap \rho_{\eta_2} \triangleleft \rho_{\eta_1} \cap \rho_{\eta_0}$ and
$\rho_{\eta_1} (\ell g(\rho_{\eta_1} \cap \rho_{\eta_2})) =
\eta_1(\ell g(\eta_1 \cap \eta_2))$.
\endroster
\enddefinition   
\bn
Note
\demo{\stag{0.5} Fact}  1) If $\lambda = \lambda^{< \kappa}$ and $E$
is ${\Cal L}$(induction)-nice \ub{then} $E$ is $L_{\lambda^+,\kappa}$-nice.
\enddemo
\newpage

\head {\S1} \endhead  \resetall \sectno=1
\bigskip

We here continue \cite[\S2]{Sh:664}, the theorem and most proofs can be
read without it.  The claims below generalize \cite{Sh:273}.
\proclaim{\stag{1.1} Claim}  Assume
\mn
$\boxtimes_1(a) \quad \lambda = \lambda^{< \lambda}$ and 
$\lambda \ge \beth_\omega$ or just
$(Dl)_\lambda$ (see \scite{0.3})
\smallskip

$(b) \quad E$ is a nice 2-place relation on $\lambda$
\smallskip

$(c)(\alpha) \quad  E$ is a equivalence relation on ${}^\lambda 2$
\smallskip

\hskip17pt $(\beta) \quad$   if $\eta,\nu \in {}^\lambda 2$ and
$(\exists!\alpha < \lambda)(\eta(\alpha) \ne \nu(\alpha))$ then $\neg
\eta E \nu$.
\mn
\ub{Then} $E$ has $2^\lambda$ equivalence classes, moreover a perfect
set of pairwise non $E$-equivalent members of ${}^\lambda 2$.
\endproclaim
\bigskip

\demo{Proof}  Note that
\mr
\item "{$\circledast$}"  of $\Bbb P$ is a $\lambda$-complete forcing
(or just $\lambda$-strategically complete) then $\Vdash_{\Bbb P}$
``clauses (c), $(\alpha), (\beta)$ are still true".  \nl
So we can apply \scite{1.1A} below.  \hfill$\square_{\scite{1.1}}$
\endroster
\enddemo
\bn
A relative is
\proclaim{\stag{1.1A} Claim}  Assume
\mn
$\boxtimes_2(a),(c) \quad$   as in $\boxtimes_1$
\smallskip

$(b)$ \hskip23pt   is $\pi^1_1[\lambda]$ 2-place relation on
${}^\lambda 2$, say defined by $(\exists z)\varphi(x,y,z,\bar a)$ \nl

\hskip40pt see Definition \scite{0.2}  
\sn

$(c)^+=(c)^+_{\text{Cohen}} \quad$   if $\Bbb P = ({}^{\lambda >}
2,\triangleleft)$, i.e. $\lambda$-Cohen then in $V^{\Bbb P}$ \nl

\hskip80pt clauses (c) still hold.
\endproclaim
\bigskip

\demo{Proof}  \ub{Stage A}:  Let $(\eta_0,\eta_1) \in {}^\lambda 2
\times {}^\lambda 2$ be generic over $\bold V$ for the forcing $\Bbb Q =
({}^{\lambda>} 2) \times ({}^{\lambda >}2)$.  Now do we have $\bold
V[\eta_0,\eta_1] \models ``\eta_0 E \eta_1"$?  If so, then for some
$(p_0,p_1) \in ({}^{\lambda >} 2) \times ({}^{\lambda >} 2)$ we have
$(p_0,p_1) \Vdash_{\Bbb Q} ``{\underset\tilde {}\to \eta_0} 
E {\underset\tilde {}\to \eta_1}"$, let
$\alpha < \lambda$ be $> \ell g(p_0),\ell g(p_1)$ and by clause
$(c)^+ (\beta)$ in $\bold V[\eta_0,\eta_1]$
we can find $\eta'_1 \in {}^\lambda 2$ such that $\eta^1_1
\restriction \alpha = \eta_1 \restriction \alpha$, and for some $\beta
\in (\alpha,\lambda),\eta'_1 \restriction [\beta,\lambda) = \eta_1
\restriction [\beta,\lambda)$, (here $\beta = \alpha +1$ is O.K.) and
$\bold V[\eta_0,\eta_1] \models \neg \eta'_1 E \eta_1$.

So $\bold V[\eta_0,\eta_1] \models ``\neg \eta_0 E \eta'_1"$ (again as in
$\bold V[\eta_0,\eta_1],E$ is an equivalence relation by clause
$(c)^+$ and $\bold V[\eta_0,\eta_1] \models ``\eta_0 E \eta_1"$).  
Also $(\eta_0,\eta'_1)$ is generic over $\bold V$ for
$({}^{\lambda >} 2) \times ({}^{\lambda >}2)$ with $(p_0,p_1)$ in the
generic set and $\bold V[\eta_0,\eta_1] = \bold V[\eta_0,\eta'_1]$ so
we get a contradiction to $(p_0,p_1) \Vdash ``\neg {\underset\tilde
{}\to \eta_0} \in {\underset\tilde {}\to \eta_1}"$.  Hence
\mr
\item "{$\circledast_1$}"  $\Vdash_{({}^{\lambda>} 2) \times ({}^{\lambda
>} 2)} ``\neg({\underset\tilde {}\to \eta_0} E
{\underset\tilde {}\to \eta_1})"$.
\endroster
\bn
\ub{Stage B}:

Let $\chi$ be large enough and let $N \prec ({\Cal H}(\chi),\in)$ be
such that $\|N\| = \lambda,N^{< \lambda} \subseteq N$ and the
definition of $E$ belongs to $N$.  Note that
\mr
\item "{$\circledast_2$}"  if $(\eta_0,\eta_1) \in ({}^\lambda 2)
\times ({}^\lambda 2)$ (and is in $\bold V$) and $N[\eta_0,\eta_1]
\models ``\neg \eta_0 E \eta_1"$, \ub{then} $\neg \eta_0 E
\eta_1$. \nl
[Why?  As $E$ is $\pi^1_1$, in $N[\eta_0,\eta_1]$ there is a witness
$\in {}^\lambda 2$ for failure, and it also witnesses in $\bold V$
that $\neg \eta_0 E \eta_1$.]
\endroster
\enddemo
\bn
Clearly to finish proving \scite{1.1}, it suffices to prove
\proclaim{\stag{1.2} Subclaim}   1) Assume $\lambda = \lambda^{< \lambda}$
and $(Dl)_\lambda$.

If ${\Cal H}(\lambda) \subseteq N,N^{<\lambda} \subseteq N,\|N\| =
\lambda$ and $N \models ZFC^-$, \ub{then} there is a perfect ${\Cal Q}
\subseteq {}^{\lambda >}2$ such that for any $\eta_0 \ne \eta_1$ from
${\Cal Q}$ the pair $(\eta_0,\eta_1)$ is generic over $N$ for
$[({}^\lambda 2) \times ({}^{\lambda >} 2)]^N$. \nl
2) Assume $\lambda = \lambda^{< \lambda}$ and there is a tree with
$\lambda$ levels each of cardinality $< \lambda$ and $2^\lambda$ \,
$\lambda$-branches.  \ub{Then} for some $X \subseteq {}^\lambda 2,|X|
= 2^\lambda$ and $\eta_0 \ne \eta_1 \in X \Rightarrow (\eta_0,\eta_1)$
generic over $N$ for $({}^{\lambda >} 2) \times ({}^{\lambda >} 2)$.
\endproclaim
\bigskip

\demo{Proof}  1) Let $\langle {\Cal P}_\alpha:\alpha < \lambda \rangle$
be such that ${\Cal P}_\alpha \subseteq {\Cal P}(\alpha),|{\Cal
P}_\alpha| < \lambda$, and for every $X \subseteq \lambda$ the set
$\{\alpha:X \cap \alpha \in {\Cal P}_\alpha\}$ is stationary.  So by
coding we can find ${\Cal P}'_\alpha \subseteq
\{(\eta_0,\eta_1):\eta_0,\eta_1 \in {}^\alpha 2\}$ of cardinality $<
\lambda$ such that for every $\eta_0,\eta_1 \in {}^\lambda 2$ the set
$\{\alpha < \lambda:(\eta_0 \restriction \alpha,\eta_1 \restriction
\alpha) \in {\Cal P}'_\alpha\}$ is stationary.  Lastly, let $\langle {\Cal
I}_\alpha:\alpha < \lambda \rangle$ list the dense open subsets of
$({}^{\lambda >} 2) \times ({}^{\lambda >} 2)$ which belongs to $N$.
Now we define by induction on $\alpha < \lambda,\langle \rho_\eta:\eta
\in {}^\alpha 2 \rangle$ such that:
\mr
\item "{$(a)$}"  $\rho_\eta \in {}^{\lambda >}2$
\sn
\item "{$(b)$}"  $\beta < \ell g(\eta) \Rightarrow \rho_{\eta
\restriction \beta} \triangleleft \rho_\eta$
\sn
\item "{$(c)$}"  $\rho_\eta \char 94 \langle \ell \rangle
\triangleleft \rho_{\eta \char 94 \langle \ell \rangle}$
\sn
\item "{$(d)$}"  if $(\eta_0,\eta_1) \in {\Cal P}'_\alpha,\ell_0 <
2,\ell_1 < 2$ then
$(\rho_{\eta_0 \char 94 \langle \ell_0 \rangle},
\rho_{\eta_1 \char 94 \langle \ell_1 \rangle}) \in \dbca_{\beta \le
\alpha} {\Cal I}_\beta$.
\ermn
There is no problem to carry the definition (using 
$|{\Cal P}'_\alpha| < \lambda = \text{ cf}(\lambda)$) and $\{\dbcu_{\alpha <
\lambda} \rho_{\eta \restriction \alpha}:\eta\in {}^\lambda 2\}$ is a
perfect set as required. \nl
2) Easier as we can assume that such a tree belongs to $N$.  
\hfill$\square_{\scite{1.2}},\square_{\scite{1.1A}}$
\enddemo
\bigskip

\proclaim{\stag{1.3} Claim}  1) In \scite{1.1}, \scite{1.1A} we can weaken clause
$(\beta)$ (in $(c), (c)^+$, call it $(c)^-,(c)^\pm$ respectively) to
\mr
\item "{$(\beta)^-$}"  if $\eta \in {}^\lambda 2,\alpha < \lambda$
then for some $\beta \in (\alpha,\lambda)$ and $\rho \in
{}^{[\alpha,\beta)} 2$ we have \nl
$\neg \eta E((\eta \restriction \alpha)
\char 94 \rho \char 94 \eta \restriction [\beta,\lambda))$.
\ermn
2) In \scite{1.1}, \scite{1.1A} and in \scite{1.3}(1) for any $\varepsilon^* \le
\lambda$ we can replace $E$ by $\langle E_\varepsilon:\varepsilon <
\varepsilon^* \rangle$, each $E_\varepsilon$ satisfying clauses (b)
and $(c),(c)^+,(c)^-,(c)^\pm$ there respectively and in the conclusion:
\mr
\item "{$(*)$}"  there is a $\lambda$-perfect ${\Cal Q}$ such that
{\roster
\itemitem{ $(\alpha)$ }  ${\Cal Q} = \langle \rho_\eta:\eta \in
{}^\lambda 2 \rangle$ and
\sn
\itemitem{ $(\beta)$ }  if $\eta_1 \ne \eta_2$ are from ${}^\lambda 2$
then $\rho_{\eta_1} \ne \rho_{\eta_2}$ and $\varepsilon < \varepsilon^*
\Rightarrow \neg(\rho_{\eta_1} E_\beta \rho_{\eta_2})$
\sn
\itemitem{ $(\gamma)$ }   for $\eta \in {}^\lambda 2$ the set $\{\ell
g(\rho_\eta \cap \rho_\nu):\nu \in {}^\lambda 2 \backslash \{\eta\}\}$
is a closed unbounded subset of $\lambda$.
\endroster}
\ermn
3) In \scite{1.1A}, \scite{1.3}(1),(2) we can weaken
$(c)^+$ or $(c)^\pm$ to
\mr
\item "{$(*)$}"  for a stationary set of $N \in [{\Cal
H}(\lambda^+)]^\lambda$ there is (in $\bold V$) \, $\eta \in
{}^\lambda 2$ which is Cohen over $N$ such that $\pi^1_1[\lambda]$
sentences are absolute from $N[\eta]$ to $\bold V$ (for
$\Sigma^1_1[\lambda]$ sentences this is necessarily true and clause
(c) (or (c)$^-$) holds.
\endroster
\endproclaim
\bigskip

\demo{Proof}  The same as the proof of \scite{1.1}.
\enddemo
\bn
Now we would like not to restrict ourselves to $\pi^1_1[\lambda]$
equivalence relations.
\proclaim{\stag{1.4} Claim}  1) Assume
\mr
\item "{$(a)$}"  $\lambda = \lambda^{< \lambda},\mu \le 2^\lambda$
\sn
\item "{$(b)$}"  $E$ is a $\pi^1_2[\lambda]$ 2-place relation on
${}^\lambda 2$, say definable by
$(\forall z_1)(\exists z_2)\varphi(x,y,z_1,z_2,a)$
{\roster
\itemitem{ $(c)(\alpha)$ }  $E$ is an equivalence relation on ${}^\lambda
2$
\sn
\itemitem{$(\beta)$ }  if $\eta,\nu \in {}^\lambda 2$ and $(\exists!
\alpha < \lambda)(\eta(\alpha)) \ne \nu(\alpha))$ then $\neg(\eta E
\nu)$
\endroster}
\item "{$(c)^+$}"  if $\eta \in {}^\lambda 2$ is generic over $\bold
V$ for $({}^{\lambda >} 2,\triangleleft)$, i.e. is Cohen \ub{then} in
$\bold V[\eta]$, clause (c) still
holds
\nl
(note that for $\rho_1,\rho_2 \in ({}^\lambda 2)^{\bold V}$ anyhow
$\bold V \models ``\rho_1 E \rho_2" \Leftrightarrow \bold V[\eta] \models
``\rho_1 E \rho_2"$)
\sn
\item "{$(d)$}"  for every $A \subseteq \lambda$ and $\chi > 2^\lambda$
there are $N,\langle \rho_\varepsilon:\varepsilon < \mu \rangle$ such
that
{\roster
\itemitem{ $(i)$ }  $N \prec ({\Cal H}(\chi),\in),N^{< \lambda} 
\subseteq N,\|N\| = \lambda,A \in N$
\sn
\itemitem{ $(ii)$ }  $\rho_\varepsilon \in {}^\lambda 2$ and $[\varepsilon <
\zeta \Rightarrow \rho_\varepsilon \ne \rho_\varepsilon]$
\sn
\itemitem{ $(iii)$ }  for $\varepsilon \ne
\zeta,(\rho_\varepsilon,\rho_\zeta)$ is generic over $N$ for
$({}^{\lambda >} 2 \times {}^{\lambda >} 2)$
\sn
\itemitem{ $(iv)$ }  $\pi^1_1[\lambda]$ formulas are preserved from
$N[\rho_\varepsilon,\rho_\zeta]$ to $\bold V$ for $\varepsilon < \zeta
< \mu$.
\endroster}
\ermn
\ub{Then} $E$ has $\ge \mu$ equivalence classes. \nl
2) We can replace $\ge \mu$ by ``perfect" in the conclusion if in (c),
$\{\rho_\varepsilon:\varepsilon < \mu\} \subseteq {}^\lambda 2$ is
perfect [see \scite{0.4}]. \nl
3) We can replace ${}^{\lambda >} 2$ by ${\Cal T} \subseteq
{}^{\lambda >} 2$ subtree such that forcing with ${\Cal T}$ adds no
bounded subset to $\lambda$.
\endproclaim
\bigskip

\demo{Proof}  By \cite[2.2t]{Sh:664}.
\enddemo
\bigskip

\definition{\stag{1.5} Definition}  Clause (d) of \scite{1.4} is called
``$\lambda$ is $[\lambda,\mu)$-w.c.a." (as in \cite[2.1t]{Sh:664}'s notation).
\enddefinition
\bigskip

\proclaim{\stag{1.6} Claim}  We can strengthen \scite{1.4} just as
\scite{1.3} strenghthens \scite{1.1}.
\endproclaim
\bn
We may wonder when does clause (d) of \scite{1.4} holds.
\proclaim{\stag{1.7} Claim}  1) Assume
\mr
\widestnumber\item{$(iii)$}
\item "{$(i)$}"  $\lambda = \lambda^{< \lambda}$ in $\bold V$
\sn
\item "{$(ii)$}"  $\Bbb P$ is a forcing notion
\sn
\item "{$(iii)$}"  $\langle {\underset\tilde {}\to
\eta_\varepsilon}:\varepsilon < \mu \rangle$ is a sequence of $\Bbb
P$-names,
\sn
\item "{$(iv)$}"  $\Vdash_{\Bbb P} ``{\underset\tilde {}\to
\eta_\zeta} \ne {\underset\tilde {}\to \eta_\varepsilon} \in
{}^\lambda 2$ for $\varepsilon < \zeta < \mu"$
\sn
\item "{$(v)$}"  if $A \subseteq \lambda,p \in \Bbb P,\chi$ large
enough then there are $N \prec ({\Cal H}(\chi),\in),\|N\| =
\lambda,N^{< \lambda} \subseteq N,\{A,p\} \in N$ and $q$ such
that $p \le q \in \Bbb P,q$ is $(N,\Bbb P)$-generic, $q \Vdash ``({}^{\lambda >}
2)^{V^{\Bbb P}} \subseteq N[{\underset\tilde {}\to G_{\Bbb P}}]"$ and
$\Bbb P' \lessdot \Bbb P$ such that $q \Vdash_{\Bbb P} ``\text{for
some } u \in [\mu]^\mu$, for every $\varepsilon \ne \zeta$ from 
$u,({\underset\tilde {}\to \eta_\varepsilon},{\underset\tilde {}\to
\eta_\zeta})$ is generic over $N[{\underset\tilde {}\to G_{{\Bbb
P}'}}]$ for $({}^{\lambda >} 2 \times {}^{\lambda >} 2)^{V^{\Bbb P}}$
and the forcing $\Bbb P/(\Bbb P + {\underset\tilde {}\to
\eta_\varepsilon} + {\underset\tilde {}\to \eta_\zeta})$ is
$\lambda$-complete (or at least $\lambda$-strategically complete).
\ermn
\ub{Then} $\lambda$ is $(\lambda,\mu)$-c.w.a. (see \scite{1.5}) in the universe $\bold
V^{\Bbb P}$.
\endproclaim
\bigskip

\demo{Proof}  Straightforward.  
\enddemo
\newpage

\head {\S2 Singulars of uncountable cofinality} \endhead  \resetall \sectno=2
\bigskip

In this section we show that usually generalization of \scite{0.1}
probably fail badly for ${}^{\text{cf}(\lambda)} \lambda,\lambda$
singular of uncountable cofinality.
\bigskip

\proclaim{\stag{2.0} Claim}  1) For every $\lambda$ strong limit
singular we can find a very nice equivalence relation $E$ on
${}^\lambda 2$ such that:
\mr
\item "{$(a)$}"  $E$ has exactly two equivalence classes
\sn
\item "{$(b)$}"  if $\eta,\nu \in {}^\lambda 2$ and $(\exists
!i)(\eta(i) \ne \nu(i))$, \ub{then} $\neg(\eta E \nu)$.
\ermn
2) For every finite abelian group $H$ and $a_\alpha \in H$ for $\alpha
< \lambda$, we can find a very nice equivalence relation $E$ on
${}^\lambda 2$ with $|H|$-equivalence classes such that the
$E$-equivalence classes can be listed as $\langle X_b:b \in H \rangle$
and
\mr
\item "{$(c)$}"  if $\eta_0,\eta_1 \in {}^\lambda 2,\eta_0(i) =
0,\eta_1(i) = 1,(\forall j)[j \ne i \Rightarrow \eta(j) =
\nu(j)]$ and $\eta_\ell \in X_{b_\ell}$ then $H \models ``b_1 - b_0
= a_i"$.
\endroster
\endproclaim
\bigskip

\remark{\stag{d.5} Remark}  1) In part (2) we can look at ${}^\lambda H$ getting a
similar result.
\endremark
\bigskip

\demo{Proof}  Use \cite[\S2]{Sh:664}.
\enddemo
\bigskip

\proclaim{\stag{2.1} Claim}  Assume
\mr
\item "{$(a)$}"  $\lambda > \kappa = \text{ cf}(\lambda) > \aleph_0$
\sn
\item "{$(b)$}"  $2^\kappa + \lambda^{< \kappa} = \lambda$.
\ermn
\ub{Then} there is $E$ such that
\mr
\item "{$(\alpha)$}"  $E$ is an equivalence relation on ${}^\kappa
\lambda$
\sn
\item "{$(\beta)$}"  $E$ is very nice (see Definition \scite{0.2})
\sn
\item "{$(\gamma)$}"  if $\eta_1,\eta_2 \in {}^\kappa \lambda$ and
$(\forall^* i < \kappa)(\eta_1(i) = \eta_2(i))$ then $\eta_1 E \eta_2
\Leftrightarrow \eta_1 = \eta_2$
\sn
\item "{$(\delta)$}"  $E$ has exactly $\lambda$ equivalence classes.
\endroster
\endproclaim 
\bigskip

\demo{\stag{2.1A} Observation}  In \scite{2.1}, and in the rest of this
section: (of course, we have to translate the results we leave it as an
exercise to the reader). \nl
1) We can restrict ourselves to $\dsize \prod_{i <
\kappa} \lambda_i,\lambda_i < \lambda = \dsize \sum_{j < \kappa}
\lambda_j$, see the proof. \nl
2) We can consider ${}^\kappa \lambda$ as a subset of ${}^\lambda 2$,
in fact a very nice one:

we identify $\eta \in {}^\kappa \lambda$ with $\nu_\eta \in {}^\kappa
2$ when $\nu_\eta(i) = 1 \Leftrightarrow i \in
\{\text{pr}(\zeta,\eta(\zeta)):\zeta < \kappa\}$ for any choice of a
pairing function pr, in fact, any one to one pr:$\kappa \times \lambda
\rightarrow \lambda$ is O.K. \nl
3) If $\lambda$ is strong limit we can identify $\dsize
\prod_{i < \kappa} \lambda_i$ with ${}^\lambda 2$ as follows:
\wilog \, $\lambda_i = 2^{\mu_i}$ with $\mu_i$ increasing, let
$\langle g^i_\varepsilon:\varepsilon < {}^{\mu_i} 2 \rangle$ list the
functions from $[\dbcu_{j < i} \mu_j,\mu_i)$ to $\{0,1\}$ and we
identify $\eta \in \dsize \prod_{i < \kappa} \lambda_i$ with $\dbcu_{i
< \kappa} g^i_{\eta(i)} \in {}^\lambda 2$. \nl
4) We can identify any $\dsize \prod_{i < \kappa} \lambda_i$ with
${}^\kappa \lambda$ when $\lambda_i \le \lambda = \dsize \sum_{i <
\kappa} \lambda_i = \text{ lim sup}\langle \lambda_i:i < \kappa
\rangle$ by identifying $\eta \in \dsize \prod_{i < \kappa} \lambda_i$
with $\nu_\eta \in {}^\kappa \lambda$ such that $1 +
\nu_\eta(\varepsilon) = \eta(\zeta)$ when $\eta(\zeta) > 0,\varepsilon
= \text{ otp}\{\xi < \zeta:\eta(\xi) > 0\}$.
\enddemo
\bigskip

\demo{Proof}  We choose $\bar \lambda = \langle \lambda_i:i < \kappa
\rangle$, nondecreasing $i < j \Rightarrow \lambda_i \le \lambda_j$
with limit $\lambda$, (e.g. $\lambda_i = \lambda$) let $\mu_j = \dsize \prod_{i < j}
\lambda_i$ so $\mu_j \le \lambda$ and let $\bar f^i = \langle
f^i_\alpha:\alpha < \mu_i \rangle$ list $\dsize \prod_{j < i}
\lambda_j$ or just a set of representatives of $\dsize \prod_{j < i}
\lambda_j/J^{bd}_i$.
\nl
For every $\eta \in {}^\kappa \lambda$ let
\mr
\item "{$(a)$}"  for limit $i < \kappa$ let $\alpha_i(\eta) = \text{ Min}\{\alpha:\eta
\restriction i = f^i_\alpha \text{ mod } J^{bd}_i\}$ 
\sn
\item "{$(b)$}"  for $\varepsilon < \kappa$ let
$B_\varepsilon(\eta) = \{i:i < \kappa$ is a limit
ordinal, $\varepsilon < i$ and $f^i_{\alpha_i(\eta)}(\varepsilon) =
\eta(\varepsilon)\}$ \nl
and lastly
\sn
\item "{$(c)$}"  $A(\eta) = \{\varepsilon <
\kappa:B_\varepsilon(\eta)$ is not stationary$\}$.  
\ermn
Now we define two binary relations $E_0,E_1$ on ${}^\kappa \lambda$:
\mr
\item "{$(d)$}"  $\eta_1 E_0 \eta_2$ \ub{iff} for every $\varepsilon <
\kappa$ we have $B_\varepsilon(\eta_1) = B_\varepsilon(\eta_2)$ mod
${\Cal D}_\kappa$, where ${\Cal D}_\kappa$ is the club filter on $\kappa$
\sn
\item "{$(e)$}" $\eta_1 E_1 \eta_2$ \ub{iff} $\eta_1 E_0 \eta_2 \and
\eta_1 \restriction A(\eta_1) = \eta_2 \restriction A(\eta_2)$.
\ermn
Clearly
\mr
\item "{$(\alpha)$}"  $E_0$ is an equivalence relation on ${}^\kappa
\lambda$ with $\le 2^\kappa < \lambda$ classes
\sn
\item "{$(\beta)$}"  $E_1$ is an equivalence relation on ${}^\kappa
\lambda$, refining $E_0$
\sn
\item "{$(\gamma)$}"  $E_0,E_1$ are very nice
\sn
\item "{$(\delta)$}"  if $\eta_1,\eta_2 \in {}^\kappa \lambda$ and
$\eta_1 E_0 \eta_2$ then $A(\eta_1) = A(\eta_2)$
\sn
\item "{$(\varepsilon)$}"  for $\eta \in {}^\kappa \lambda,A(\eta_1)$
is a bounded subset of $\kappa$ \nl
[why?  otherwise let $C = \{\delta < \kappa:\delta = \sup(A(\eta_1)
\cap \delta)\}$, it is a club of $\kappa$, and for each $i \in C$ there
is $j_i < i$ such that $\eta \restriction [j_i,i) =
f^i_{\alpha_i(\eta)} \restriction [j_i,i)$, clearly $j_i$ exists by the
defintion of $\alpha_i(\eta)$.  By Fodor lemma, for some $j(*) <
\kappa$ the set $S_{j(*)} = \{i \in C:j_i = j(*)\}$ is stationary, now
choose $\varepsilon \in A(\eta) \backslash j(*)$, so clearly
$B_\varepsilon(\eta)$ includes $S_{j(*)}$ hence is a stationry subset
of $\kappa$ hence by the definition of $A(\eta)$ clearly $\varepsilon$
does not belong to $A(\eta)$, contradiction.]
\nl
So clearly
\sn
\item "{$(\zeta)$}"  $E_1$ has $\le ({}^\kappa \lambda/E_0) + \Sigma\{
\dsize \prod_{j < i} \lambda_j:i <\kappa\} \le \lambda$ equivalence
class.
\nl
Now
\item "{$(\eta)$}"  if $\eta_1,\eta_2 \in {}^\kappa \lambda$ and
$\eta_1 = \eta_2 \text{ mod } J^{bd}_\kappa$ \ub{then} for every limit
$i < \kappa$ large enough we have $\alpha_i(\eta_1) =
\alpha_i(\eta_2)$ \nl
[why?  let $i^* = \sup\{j+1:\eta_1(j) \ne \eta_2(j)\}$ so by the
assumption, if $i$ limit $\and i \in (i^*,\kappa)$ then $\eta_1
\restriction i = \eta_1 \restriction i$ mod $J^{bd}_i$ hence
$\alpha_i(\eta_1) = \alpha_i(\eta_2)$ by the definition of
$\alpha_i(-)$, which is the desired conclusion of $(\eta)$.]
\sn
\item "{$(\theta)$}"  if $\eta_1,\eta_2 \in {}^\kappa \lambda$ and
$\eta_1 = \eta_2$ mod $J^{bd}_\kappa$ \ub{then} $\eta_1 E_1 \eta_2
\Leftrightarrow \eta_1 = \eta_2$ \nl
[why?  if $\eta_1 = \eta_2$ clearly $\eta_1 E_1 \eta_2$; so assume
$\eta_1 E_1 \eta_2$ and we shall show that $\eta_1 = \eta_2$,
i.e. $\varepsilon < \kappa \Rightarrow \eta_1(\varepsilon) =
\eta_2(\varepsilon)$.  By the definition of $E_1$ we have $\eta_1 E_0
\eta_2$ hence by clause $(\delta)$ we have $A(\eta_1) = A(\eta_2)$,
call it $A$.  If $\varepsilon \in A$, by the definition of $E_1$ we have
$\eta_1 \restriction A = \eta_2 \restriction A$ hence
$\eta_1(\varepsilon) = \eta_2(\varepsilon)$.  So assume $\varepsilon \in
\kappa \backslash A$, first we can find $j^* < \kappa$ such that for every
limit $i \in (j^*,\kappa)$ we have $\alpha_i(\eta_1) =
\alpha_i(\eta_2)$, it exists by clause $(\eta)$.  Second, then
$B_\varepsilon(\eta_1),B_\varepsilon(\eta_2)$ are stationary (as
$\varepsilon \notin A(\eta_\ell)$) and equal modulo ${\Cal D}_\kappa$
(as $\eta_1 E_0 \eta_2$); so we can find $i \in B_\varepsilon(\eta_1)
\cap B_\varepsilon(\eta_2)$ which satisfy $i > j^*$.  Now $\eta_1(\varepsilon)
= f^i_{\alpha_i(\eta_1)}(\varepsilon)$ by the definition of
$B_\varepsilon(\eta_1)$ as $i \in B_\varepsilon(\eta_1)$ and
$\alpha_i(\eta_1) = \alpha_i(\eta_2)$ as $i > j^*$ and
$f^i_{\alpha_i(\eta_2)}(\varepsilon) = \eta_2(\varepsilon)$ by the
definition of $B_\varepsilon(\eta_2)$ as $i \in
B_\varepsilon(\eta_1)$; together $\eta_1(\varepsilon) =
\eta_2(\varepsilon)$.  So we have completed the proof that $\varepsilon
< \kappa \Rightarrow \eta_1(\varepsilon) = \eta_2(\varepsilon)$ thus
proving $\eta_1 = \eta_2$ as required.]
\sn
\item "{$(\iota)$}"  $E_1$ has $\ge \lambda_i$ equivalence classes for
any $i < \kappa$ \nl
[why? let $\eta^* \in {}^\kappa \lambda$ and for $\alpha < \lambda_i$
let $\eta^*_\alpha \in {}^\kappa \lambda$ be defined by
$\eta^*_\alpha(\varepsilon)$ is $\alpha$ if $\varepsilon = i$ and is
$\eta^*(\varepsilon)$ otherwise.  By clause $(\theta)$ we have $\alpha
< \beta < \lambda_i \Rightarrow \neg \eta^*_\alpha E_1 \eta^*_\beta$,
hence $|{}^\kappa \lambda/E_1| \ge \lambda_i$.]
\sn
\item "{$(\kappa)$}"  $E_1$ has exactly $\lambda$ equivalence classes
\nl
[why?  by clause $(\iota), E_1$ has $\ge \sup\{\lambda_i:i < \kappa\}
= \lambda$ equivalence classes and by clause $(\zeta)$ has $\le
\lambda$ equivalence classes.]  \hfill$\square_{\scite{2.1}}$
\endroster
\enddemo
\bigskip

\proclaim{\stag{2.2} Claim}  Assume
\mr
\item "{$(a)$}"  $\lambda > \kappa = \text{ cf}(\lambda) > \aleph_0$
\sn
\item "{$(b)$}"  $2^\kappa + \lambda^{< \kappa} = \lambda$
\sn
\item "{$(c)$}"  $\lambda \le \theta \le \lambda^\kappa$.
\ermn
\ub{Then} there is $E$ such that
\mr
\item "{$(\alpha)$}"  $E$ is an equivalence relation on ${}^\kappa
\lambda$
\sn
\item "{$(\beta)$}"  $E$ is very nice
\sn
\item "{$(\gamma)$}"  if $\eta_1,\eta_2 \in {}^\kappa \lambda$ and
$\eta_1 = \eta_2$ mod $J^{bd}_\kappa$ then $\eta_1 E \eta_2
\Leftrightarrow \eta_1 = \eta_2$
\sn
\item "{$(\delta)$}"  $E$ has $\theta$ equivalence classes.
\endroster
\endproclaim
\bigskip

\demo{Proof}  We can find a tree ${\Cal T} \subseteq {}^\kappa
\lambda$ with $\lambda$ nodes and exactly $\theta$ \,
$\kappa$-branches (\cite{Sh:262}); we can easily manage that $\eta \ne
\nu \in \text{ lim}_\kappa({\Cal T}) \Rightarrow (\exists^\kappa i <
\kappa)(\eta(i)) \ne \nu(i))$.  We proceed as in the proof of
\scite{2.1}, but in the definition of $E_1$ we add

$$
\eta_1 \in \text{ lim}_\kappa(T) \equiv \eta_2 \in \text{
lim}_\kappa({\Cal T}) \and (\eta_1 \in \text{ lim}_\kappa({\Cal T})
\rightarrow \eta_1 = \eta_2).
$$
${{}}$  \hfill$\square_{\scite{2.2}}$
\enddemo
\bigskip

\proclaim{\stag{2.3} Claim}  In Claim \scite{2.1} we replace clause
$(\gamma)$ by
\mr
\item "{$(\gamma)_1$}"  for every $\eta^* \in {}^\kappa \lambda$,
the set $\{\eta \in {}^\kappa \lambda:\eta = \eta^* \text{ mod }
J^{bd}_\kappa\}$ is a set of representatives for the family of
$E$-equivalence classes.
\endroster
\endproclaim
\bigskip

\demo{Proof}  Let $\bar \lambda$ be as there.

Let $K_i$ be an additive group, with set of elements $\lambda_i$ and
zero $0_{K_i}$.  Let $<^*$ be a well ordering of ${}^\kappa({\Cal P}(\kappa))$.

For every $\eta \in \dsize \prod_{i < \kappa} \lambda_i$ let

$$
\Xi_\eta = \{\langle B_\varepsilon(\nu):\varepsilon < \kappa
\rangle:\nu \in \dsize \prod_{i < \kappa} \lambda_i \text{ and } 
\nu = \eta \text{ mod } J^{bd}_\kappa\}
$$
\mn
so it is a nonempty subset of ${}^\kappa({\Cal P}(\kappa))$ and let
$\bar B^*_\eta$ be its $<^*$-first member and let $\Theta_\eta = \{\nu \in
\dsize \prod_{i < \kappa} \lambda_i:B_\varepsilon(\nu) = B^*_{\eta,\varepsilon}$ for every
$\varepsilon < \kappa$ and $\nu = \eta$ mod $J^{\text{bd}}_\kappa\}$.
\sn
Now note
\mr
\item "{$(*)_0$}"  $\Theta_\eta \ne \emptyset$ so $\bar B^*_\eta$ is
well defined for $\eta \in \dsize \prod_{i < \kappa} \lambda_i$
\sn
\item "{$(*)_1$}"  if $\nu \in \Theta_\eta$ \ub{then} for every limit $i <
\kappa$ large enough we have $\alpha_i(\nu) = \alpha_i(\eta)$
\sn
\item "{$(*)_2$}"  if $\nu_1,\nu_2 \in \Theta$ and $\varepsilon <
\kappa$, \ub{then} for every limit $i$ large enough we have:
$f^i_{\alpha_i(\nu_1)}(\varepsilon) =
f^i_{\alpha_i(\nu_2)}(\varepsilon)$.
\ermn
Now for $\eta \in \dsize \prod_{i < \kappa} \lambda_i$ we define
$\rho_\eta \in \dsize \prod_{i < \kappa} \lambda_i$ by

$$
\align
\rho_\eta(\varepsilon) \text{ is}: &f^i_{\alpha_i(\eta)}(\varepsilon)
\text{ for every } i \in B^*_\varepsilon \\
  &\text{ large enough if } B^*_{\eta,\varepsilon} \text{ is stationary}
\endalign
$$

$$
0_{K_i} \text{ if } B^*_\varepsilon \text{ is not stationary}.
$$
\mn
It is easy to see that
\mr
\item "{$(*)_3$}"  if $\nu \in \kappa_\lambda$ then $\rho_\eta(\varepsilon) =
\eta(\varepsilon)$ for every $\varepsilon < \kappa$ large enough hence
\sn
\item "{$(*)_4$}"  $\rho_\eta = \eta$ mod $J^{bd}_\kappa$
\sn
\item "{$(*)_5$}"  if $\eta_1,\eta_2 \in \dsize \prod_{i < \kappa}
\lambda_i$ and $\eta_1 = \eta_2 \text{ mod } J^{bd}_\kappa$ \ub{then}
$\rho_{\eta_1} = \rho_{\eta_2}$.
\ermn
Lastly, we define the equivalence relation $E$:

for $\eta_1,\eta_2 \in \dsize \prod_{i < \kappa} \lambda_i$ we define:

$$
\eta_1 E \eta_2 \text{ \ub{iff} (for every } i < \kappa \text{ we
have } K_i \models ``\eta_1(i)-\rho_{\eta_1}(i) =
\eta_2(i)-\rho_{\eta_2}(i)").
$$ 
\mn
Now check.  \hfill$\square_{\scite{2.3}}$
\enddemo
\bn
We may like to weaken the cardinal arithmetic assumptions.
\proclaim{\stag{2.4} Claim}  Assume
\mr
\item "{$(a)$}"  $\lambda > \kappa = \text{ cf}(\lambda) > \aleph_0$
\sn
\item "{$(b)$}"  $\kappa^{\aleph_0} < \lambda = \lambda^{\aleph_0}$.
\ermn
\ub{Then} the results \scite{2.1}, \scite{2.3} and \scite{2.2} (if
there is a tree ${\Cal T}$ with $\lambda$ nodes and $\kappa$-branches)
holds \ub{if} we replace the ideal $J^{bd}_\kappa$ by the ideal
$[\kappa]^{< \aleph_0}$.
\endproclaim
\bigskip

\demo{Proof}  Without loss of generality $\lambda_i >
\kappa^{\aleph_0},\langle \lambda_i:i < \kappa \rangle$ as in \scite{2.1}.
Let $\langle D_i:i < \kappa^{\aleph_0} \rangle$ list the
subsets of $\kappa$ of order type $\omega$ and let $\bar f^i = \langle
f^i_\alpha:\alpha < \dsize \prod_{j \in D_i} \lambda_j \rangle$ list
$\dsize \prod_{j \in D_i} \lambda_j$.  For $\eta \in \dsize
\prod_{\varepsilon < \kappa} \lambda_\varepsilon$ let $\alpha_i(\eta)
= \text{ Min}\{\alpha:\eta \restriction D_i = f^i_\alpha \text{ mod
} J^{bd}_{D_i}\}$.  With those choices the proofs are similar, we
write below the proof of the parallel of \scite{2.1}.

Define a two place relation $E^*$ on $\dsize \prod_{\varepsilon <
\kappa} \lambda_\varepsilon:\eta E \nu$ iff $(\forall i <
\kappa^{\aleph_0})(\alpha_i(\eta) = \alpha_i(\nu))$. Clearly $E^*$ is
a very nice equivalence relation and $\eta E^* \nu$ iff $(\eta,\nu
\in \dsize \prod_{\varepsilon < \kappa} \lambda_\varepsilon$ and
$\{\varepsilon:\eta(\varepsilon) \ne \nu(\varepsilon)\}$ is finite).
For $\eta \in \dsize \prod_{\varepsilon < \kappa} \lambda_\varepsilon$
let $A_{\eta,\varepsilon} = \{f^i_{\alpha_i(\eta)}(\varepsilon):i <
\kappa^{\aleph_0}$ and $\varepsilon \in D_i\}$ and let
$\gamma_{\eta,\varepsilon}$ be $0$ if $\eta(\varepsilon) \in
A_{\eta,\varepsilon}$ and $\eta(\varepsilon)-\sup(A_{\eta,\varepsilon}
\cap \eta(\varepsilon))$ otherwise stipulating $\alpha - \beta = 0$ if
$\alpha \le \beta$.  Lastly, $\eta E \nu$ iff
$(\forall \varepsilon)(\gamma_{\eta,\varepsilon} \equiv
\gamma_{\nu,\varepsilon})$.
\hfill$\square_{\scite{2.4}}$
\enddemo
\bigskip

\proclaim{\stag{2.5} Claim}  1) If 
$2^\kappa < \lambda = \lambda^{\aleph_0},\aleph_0 < \kappa =
\text{ cf}(\lambda) < \lambda$, \ub{then} we can find $E$ as in
\scite{2.1}$(\alpha),(\beta)$ (but not $(\gamma)$) and
\mr
\item "{$(\gamma)^*$}"  if $\eta \in {}^\kappa \lambda$ and $i <
\kappa$ then $X_{\alpha,i} = \{\nu \in {}^\kappa \lambda:(\forall j)(j
< \kappa \and j \ne i \rightarrow \nu(j) = \eta(j)\}$ is a set of
representations for $E$.
\ermn
2)  If $2^\kappa \le \lambda = \lambda^{\aleph_0},\aleph_0 < \kappa =
\text{ cf}(\lambda) < \lambda,1 \le \theta \le \lambda$, \ub{then} we
can find $E$ as in \scite{2.1}$(\alpha),(\beta)$ and
\mr
\item "{$(\gamma)^*$}"  if $\eta \in {}^\kappa \lambda$ and $i <
\kappa$ then $X_{\alpha,i}$ contains a set of representations
\sn
\item "{$(\delta)^*$}"   $E$ has $\theta$ equivalence classes.
\endroster
\endproclaim
\bigskip

\demo{Proof}  1)  We let $K$ be an abelian group with universe
$\lambda$ and in the proof of \scite{2.3},
define $\rho_\eta$ such that $\eta = \rho_\eta$ mod
$[\kappa]^{< \aleph_0}$ and let $a_\eta = \{i < \kappa:\eta(i) \ne
\rho_\eta(i)\} \in [\kappa]^{< \aleph_0}$ and define $E$ by: $\eta E
\nu$ iff $K \models \dsize \sum_{i \in a_\eta} \eta(i) -
\rho_\eta(i) = \dsize \sum_{i \in a_\nu} \nu(i) - \rho_\nu(i)$. \nl
2) Similar but we use equality in $K/K_1,K_1$ a subgroup of $K$ such that
$(K:K_1) = \theta$.  \hfill$\square_{\scite{3.5}}$
\enddemo
\bigskip 

\remark{\stag{2.6} Concluding Remark}  Instead $\langle J^{bd}_{D_i}:i
< \kappa^{\aleph_0} \rangle$ we can use $\langle (D_i,J_i):i < i^{\bar
\lambda} \rangle,D_i \subseteq \kappa,J_i$ an ideal on $D_i$ such that
$|\dsize \prod_{\varepsilon \in D_i} \lambda_\varepsilon/J_i| \le
\lambda,I = \{D \subseteq \kappa$: for every $i < i^*$ we have $D \cap
D_i \in J_i\}$ is included in $J^{bd}_\kappa$. \nl
Have not pursued this.
\endremark
\newpage

\head {\S3 Countable cofinality: positive results} \endhead  \resetall \sectno=3
\bigskip

We first phrase sufficient conditions which related to large
cardinals.  Then we prove that they suffice.  The proof of \scite{3.1}
is later in this section.
\proclaim{\stag{3.1} Lemma} Assume
\mr
\item "{$(a)$}"  $\lambda$ is strong limit of cofinality $\aleph_0$
\sn
\item "{$(b)$}"  $\lambda$ is a limit of measurables, or just 
\sn
\item "{$(b)^-$}"  for every $\theta < \lambda$ for some $\mu,\chi$
satisfying $\theta \le \mu \le \chi < \lambda$, 
there is a $(\chi,\mu,\theta)$-witness (see Definition \scite{3.2} below) 
\sn
\item "{$(c)$}"  $E$ is a nice equivalence relation on ${}^\omega
\lambda$ (or has enough absoluteness, as proved in \scite{3.11})
\sn
\item "{$(d)$}"  if $\eta,\nu \in {}^\omega \lambda$ and
$(\exists!n)(\eta(n) \ne \nu(n))$ then $\neg(\eta E \nu)$.
\ermn
\ub{Then} $E$ has $2^\lambda$ equivalence classes, moreover if
$\lambda_n < \lambda_{n+1} < \lambda = \dsize \Sigma_{n < \omega} \,
\lambda_n$ \ub{then} there is a subtree of
${}^{\omega >} \lambda$ isomorphic to $\dbcu_m \dsize \prod_{n < m}
\lambda_n$, with the $\omega$-branches pairwise non $E$-equivalent
(even somewhat more, see \scite{3.14}). 
\endproclaim
\bigskip

\definition{\stag{3.2} Definition}  1) We say $(\Bbb Q,I^{s_1,s_2}_2)$ is a
$(\lambda,\mu,\theta)$-witness if ($\lambda \ge \mu \ge \theta$ and):
\mr
\item "{$(a)$}"  $\Bbb Q$ is a $\theta$-complete forcing notion
\sn
\item "{$(b)$}"  $s_1$ is a function from $\Bbb Q$ to ${\Cal P}(\lambda)
\backslash \{\emptyset\}$
\sn
\item "{$(c)$}"  $s_2$ is a function from $\Bbb Q$ to $\{A:A \subseteq
\{(\alpha,\beta):\alpha < \beta < \lambda\}\}$
\sn
\item "{$(d)$}"  if $\Bbb Q \models p \le q$ then $s_\ell(q) \subseteq
s_\ell(p)$ for $\ell = 1,2$
\sn
\item "{$(e)$}"  $(\alpha,\beta) \in s_2(p) \Rightarrow
\{\alpha,\beta\} \subseteq s_1(p)$
\sn
\item "{$(f)$}"  for every $p \in \Bbb Q$ there is $q$ such that $p \le q
\in \Bbb Q$ and \nl
$(\forall \beta)(\exists \alpha,\gamma)[\beta \in s_1(q) \rightarrow
(\alpha,\beta) \in s_2(p) \and (\beta,\gamma) \in s_2(p)]$
\sn
\item "{$(g)$}"  if $p \in \Bbb Q$ and $A \subseteq \lambda \times
\lambda$, \ub{then} for some $q$ we have $p \le q \in \Bbb Q$ and $s_2(p)
\subseteq A \vee s_2(p) \cap A = \emptyset$
\sn
\item "{$(h)$}"  if $p \in \Bbb Q$ then for some $Y \in [\lambda]^\mu$
for every $\alpha < \beta$ from $Y$ we have $(\alpha,\beta) \in s_2(p)$. 
\ermn
2) We say $(\Bbb Q,s_1,s_2)$ is a $(\lambda,\mu,\theta,\varrho)$-witness if we
can strengthen clause $(g)$ to
\mr
\item "{$(g)^+$}"  if $f:{}^2 \lambda \rightarrow \varrho$ 
and $p \in \Bbb Q$ \ub{then} for some $q$ we have $p \le q \in
\Bbb Q$ and $f \restriction s_2(q)$ is constant.
\ermn
3) We call $(\Bbb Q,s_1,s_2)$ uniform if $\lambda =: \cup\{s_1(p):p
\in \Bbb Q\}$ is a cardinal and for every $p \in \Bbb Q$ and $\alpha < \lambda$ for some
$q$ we have $p \le q \in \Bbb Q$ and $s_1(q) \cap \alpha =
\emptyset$. \nl
4) We replace $\varrho$ by $< \varrho$ if in clause $(g)^+$, Rang$(f)$
is a subset of $\varrho$ of cardinality $< \varrho$.
\enddefinition
\bigskip

\definition{\stag{3.3} Definition}  1) We say that $(\Bbb Q,\bar s)$ is a
$(\lambda,\mu,\theta,\varrho;n)$-witness if $\lambda \ge \mu \ge
\theta,\lambda \ge \varrho$ and $\bar s = \langle s_m:m =
1,\dotsc,n \rangle$ and
\mr
\item "{$(a)$}"  $\Bbb Q$ is a $\theta$-complete forcing
\sn
\item "{$(b)$}"  $s_m$ is a function from $\Bbb Q$ to ${\Cal P}(\{\bar
\alpha:\bar \alpha = \langle \alpha_\ell:\ell < m \rangle,\alpha_\ell
< \alpha_{\ell +1}$ for $\ell < m-1\})$
\sn
\item "{$(c)$}"  if $\Bbb Q \models p \le q$ and $m \in
\{1,\dotsc,n\}$ then $s_m(p) \subseteq s_m(q)$
\sn
\item "{$(d)$}"  if $\langle \alpha_\ell:\ell < m+1 \rangle \in s_{m+1}(p)$
and $k < m+1$ then \nl
$\langle \alpha_\ell:\ell < k \rangle \char 94 \langle
\alpha_\ell:\ell = k+1,\dotsc,m \rangle \in s_m(p)$
\sn
\item "{$(e)$}"  for every $m \in \{1,\dotsc,n-1\},k<m$ and $p \in
\Bbb Q$ there is $q$ satisfying $p \le q \in \Bbb Q$ and $(\forall \bar \alpha \in
s_m(q))(\exists \bar \beta \in s_{m+1}(p))[\bar \alpha = (\bar \beta
\restriction k) \char 94 (\bar \beta \restriction [k+1,m))]$
\sn
\item "{$(f)^+$}"  if $m \in \{1,\dotsc,n\}$ and $f:{}^m \lambda
\rightarrow \varrho$ and $p \in \Bbb Q$ \ub{then} for some $q$ we have
$p \le q \in \Bbb Q$ and $f \restriction s_m(q)$ is constant
\sn
\item "{$(g)$}"  if $p \in \Bbb Q$ then for some $Y \in [\lambda]^\mu$
every increasing $\bar \alpha \in {}^n Y$ belongs to $s_n(p)$.
\ermn
2) $(\Bbb Q,\bar s)$ is a $(\lambda,\theta,\varrho;\omega)$-witness 
is defined similarly (using $s_m$ for $m \in [1,\omega)$). \nl
3) If $\varrho = 2$ we may omit it, also in Definition \scite{3.2}.
Also ``uniform" and ``$< \varrho$" means as in Definition \scite{3.2}. 
\enddefinition
\bigskip

\proclaim{\stag{3.4} Claim}  1) If $(\Bbb Q,s_1,s_2)$ is a
$(\lambda,\mu,\theta;n)$-witness and $\varrho < \theta,n < \omega$,
\ub{then} $(\Bbb Q,s_1,s_2)$ is a
$(\lambda,\mu,\theta,2^\varrho;n)$-witness. \nl
2) If ${\Cal D}$ is a normal ultrafilter on $\lambda$ so $\lambda$ is
measurable, $\Bbb Q = ({\Cal
D},\supseteq),s_1(A) = A,s_2(A) = \{(\alpha,\beta):\alpha < \beta$ are
from $A\}$, \ub{then} $(Q,s_1,s_2)$ is a uniform
$(\lambda,\lambda,\lambda)$-witness. \nl
3) If in (2), $s_m(A) = \{\bar \alpha:\bar \alpha = \langle
\alpha_\ell:\ell < m \rangle$ is increasing, $\alpha_\ell \in A\},\bar
s = \langle s_{m+1}:m < n \rangle$ and $n \le \omega$ \ub{then} $(\Bbb Q,\bar s)$ is a
$(\lambda,\lambda,\lambda,< \lambda;n)$-witness. \nl
4) If there is a $(\lambda,\mu,\theta;n)$-witness and $2^{< \theta} \le
\lambda$, \ub{then} there is such $(\Bbb Q,\bar s)$ with $|\Bbb Q| \le
2^\lambda$. \nl
5) Definition \scite{3.2}(1) is the case $n=2$ of Definition
\scite{3.3}(1).
\endproclaim
\bigskip

\demo{Proof}  Easy.
\enddemo
\bigskip

\proclaim{\stag{3.5} Claim}  1) Assume
\mr
\item "{$(a)$}"  $2 \le n < \omega,\lambda = \beth_{n-1}(\theta)^+$
\sn
\item "{$(b)$}"  $\theta$ is a compact cardinal or just a $\lambda$-compact cardinal
\sn
\item "{$(c)$}"  $\mu = \mu^{< \mu} < \theta$
\sn
\item "{$(d)$}"  $\Bbb P = \text{ Levy}(\mu,< \theta)$.
\ermn
\ub{Then} in $V^{\Bbb P}$, there is a $(\lambda,\mu,\theta;n)$-witness
$(\Bbb Q,\bar s)$. \nl
2) If there are $\lambda_n$ for $n < \omega,\lambda_n < \lambda_{n+1}$
and $\lambda_n$ is
$2^{(2^{\lambda_n})^+}$-compact, \ub{then} for some set forcing $\Bbb P$ in
$V^{\Bbb P}$ the cardinal $\lambda = \beth_\omega = \aleph_\omega$ is
dichotomically good (see Definition \scite{3.6A} below).
\endproclaim
\bigskip

\demo{Proof}  By \cite{Sh:124}.
\enddemo
\bn
Toward proving Lemma \scite{3.1} assume (till the end of this section) that
\demo{\stag{3.6} Hypothesis}  $\lambda = \dsize \sum \lambda_n,\dsize
\sum_{\ell < n} 2^{\lambda_\ell} < \theta_n \le \lambda_n$ and $(\Bbb
P_n,s_{n,1},s_{n,2})$ is a $(\lambda_n,\mu^+_n,\theta_n)$-witness and
for simplicity $\mu_n < \mu_{n+1},\lambda = \dsize \sum_n \mu_n$.
\enddemo
\bigskip

\definition{\stag{3.6A} Definition}  We call $\lambda$ dichotomically
good if there are $\lambda_n,\mu_n,\theta_n,\Bbb P_n,s_{n,1},s_{n,2}$
as in \scite{3.6}.
\enddefinition
\bigskip

\definition{\stag{3.7} Definition}  1) We define the forcing notion
$\Bbb Q_1$
\mr
\item "{$(a)$}"  $\Bbb Q_1 = \bigl\{p:p = (\eta,\bar A) = (\eta^p,\bar A^p) \text{ such
that letting } n(p) = \ell g(\eta)$ \nl
\hskip80pt $\text{we have } n^p < \omega,
\eta^p \in \dsize \prod_{\ell < n[p]} \lambda_\ell$ and \nl
\hskip80pt $\bar A^p = \langle A^p_\ell:\ell \in [n(p),\omega) \rangle \text{ and } 
A^p_\ell \in \Bbb P_\ell \bigr\}$
\sn
\item "{$(b)$}"  $p \le_{{\Bbb Q}_1} q$ iff $\eta^p \trianglelefteq
\eta^q,n(p) \le \eta(q),\ell \in [n(q),\omega) \Rightarrow \Bbb P_\ell
\models ``A^p_\ell \le A^q_\ell"$ and $n(p) \le \ell < n(q)
\Rightarrow \eta^q(\ell) \in s_1(A^p_\ell)$
\sn
\item "{$(c)$}"  We define the $\Bbb Q_1$-name $\underset\tilde {}\to
\eta$ by: $\underset\tilde {}\to \eta[G] = \cup\{\eta^p:
p \in {\underset\tilde {}\to G_{{\Bbb Q}_1}}\}$
\sn
\item "{$(d)$}"  We define
{\roster
\itemitem{$ (\alpha)$ }   $p \le^{{\Bbb Q}_1}_{pr} q$ iff $p
\le_{{\Bbb Q}_1} q \and n(p) = n(q)$
\sn
\itemitem{ $(\beta)$ }  $p \le^{{\Bbb Q}_1}_{apr} q$ iff
$p \le_{{\Bbb Q}_1} q \and \dsize \bigwedge_{\ell \ge n(q)} A^q_\ell =
A^p_\ell$
\sn
\itemitem{ $(\gamma)$ }  $p \le^{{\Bbb Q}_1}_{pr,n} q$ iff $p
\le^{{\Bbb Q}_1}_{pr} q$ and $\bar A^p \restriction [n(p),n(q)) = \bar
A^q \restriction [n(p),n(q))$.
\endroster}
\ermn
2) We define the forcing notion $\Bbb Q_2$ by: 
\mr
\item "{$(a)$}" $\Bbb Q_2 = \bigl\{p:p = (\eta_0,\eta_1,\bar A) =
(\eta^p_0,\eta^p_1,\bar A^p) \text{ where for some } n(p) < \omega
\text{ we have}$:
\nl
\hskip40pt  $\eta^p_0,\eta^p_1 \in \dsize \prod_{\ell < n[p]}
\lambda_\ell,\bar A^p = \langle A^p_\ell:\ell \in [n(p),\omega)
\rangle$ and $A^p_\ell \in \Bbb P_\ell \bigr\}$
\sn
\item "{$(b)$}"  $p \le_{{\Bbb Q}_2} q$ iff
{\roster
\itemitem{ $(i)$ }   $n(p) \le n(q)$ 
\sn
\itemitem{ $(ii)$ }  $\eta^p_\ell \trianglelefteq \eta^q_\ell$ for
$\ell=0,1$ 
\sn
\itemitem{ $(iii)$ }  $A^q_\ell \subseteq A^p_\ell$ for $\ell \in
[n(q),\omega)$ 
\sn
\itemitem{ $(iv)$ }  the pair 
$(\eta^p_0(\ell),\eta^p_1(\ell))$ is from $s_2(A^p_\ell)$ for $\ell
\in [n(p),n(q))$
\endroster}
\item "{$(c)$}"  we define
\footnote{why don't we ask $(\forall \ell <
n^p)(\eta^p_0(\ell) < \eta^p_1(\ell))$? to be able to construct the
perfect set, but, of course, $p \Vdash_{{\Bbb Q}_2} ``{\underset\tilde
{}\to \eta_0}(\ell) < {\underset\tilde {}\to \eta_1}(\ell)$ for $\ell
\in [n(p),\omega)$"} the $\Bbb Q_2$-name ${\underset\tilde {}\to
\eta_\ell}$ (for $\ell =0,1)$ by 
${\underset\tilde {}\to \eta_\ell}[G] = \cup\{\eta^p_\ell:p \in
{\underset\tilde {}\to G_{{\Bbb Q}_2}}\}$
\sn
\item "{$(d)$}"  we define
{\roster
\itemitem{ $(\alpha)$ }   $p \le^{{\Bbb Q}_2}_{pr} q$ iff $p
\le_{{\Bbb Q}_1} q \and n(p) = n(q)$ and
\sn
\itemitem{$ (\beta)$ }   $p \le^{{\Bbb Q}_2}_{apr} q$
iff $p \le_{{\Bbb Q}_2} q \and \dsize \bigwedge_{\ell \ge n(q)}
A^q_\ell = A^p_\ell$ and
\sn
\itemitem{ $(\gamma)$ }  $p \le^{{\Bbb Q}_2}_{pr,n} q$ iff
$p \le^{{\Bbb Q}_2}_{pr} q$ and $\bar A^p \restriction [n(p),n) \equiv
\bar A^q \restriction [n(p),n)$.
\endroster}
\ermn
3) If for a fix $k < \omega$, we have 
$(\Bbb P_n,\bar s^n)$ is a $(\lambda_n,\mu_n,\theta_n;k)$-witness for $n <
\omega$ \ub{then} we can define $\Bbb Q_\ell$ for $\ell = 1,\dotsc,k$ naturally. \nl
4) If $(\Bbb P_n,\bar s^n)$ is a $(\lambda_n,\mu_n,\theta_n;n)$-witness for $n <
\omega$ \ub{then} we can define $\Bbb Q = \{(\eta_0,\bar A):n < \omega,\eta(\ell) \in
{}^\ell(\lambda_\ell)\}$ is increasing, $\bar A = \langle A_\ell:\ell \in
[n,\omega),A_\ell \in {\Bbb P}_\ell\}$ with the natural order. \nl
We shall not pursue (3), (4).
\enddefinition
\bigskip

\demo{\stag{3.9} Fact}  Let $\ell \in \{1,2\}$. \nl
1) If $p \le_{{\Bbb Q}_\ell} q$ \ub{then} 
for some $q$ we have $p \le^{{\Bbb Q}_\ell}_{pr,n(q)} q \le^{{\Bbb Q}_\ell}_{apr}
r$. \nl
2) If $\bar p = \langle p_i:i < \alpha \rangle$ is $\le^{{\Bbb
Q}_\ell}_{pr}$-increasing, $\alpha < \theta_{n(p_0)}$, \ub{then}
$\bar p$ has a $\le^{{\Bbb Q}_\ell}_{pr}$-upper bound. \nl
3) If $\underset\tilde {}\to \tau$ is a ${\Bbb Q}_\ell$-name of an
ordinal, $p \in {\Bbb Q}_\ell$, \ub{then} for some $q$ and $n$ we
have:
\mr
\item "{$(a)$}"  $p \le_{pr} q$
\sn
\item "{$(b)$}"  if $q \le_{apr} r$ and $n(r) \ge n$, \ub{then} $r$ forces a
value to $\underset\tilde {}\to \tau$.
\ermn
4) In (3), if $\Vdash ``\underset\tilde {}\to \tau < \omega$ or just
$< \alpha^* < \theta_{n(p)}"$ then \wilog \, $n = n(p)$.
\enddemo
\bigskip

\demo{Proof}  Easy.
\enddemo
\bigskip

\proclaim{\stag{3.10} Claim}  1) $\underset\tilde {}\to \eta$ is
generic for ${\Bbb Q}_1$. \nl
2) The pair $({\underset\tilde {}\to \eta_0},{\underset\tilde {}\to \eta_1}), 
{\underset\tilde {}\to \eta_1}$ is generic for $\Bbb Q_2$.
\endproclaim
\bigskip

\proclaim{\stag{3.11} Claim}  Forcing by $\Bbb Q_2$ preserve ``$E$ is a
nice (see Definition \scite{0.2}(2)) equivalence relation on $\dsize
\prod_{n < \omega} \lambda_n$ satisfying clause (d) of \scite{3.1}.
\endproclaim
\bigskip

\demo{Proof}  Assume toward contradiction that $p^* \Vdash_{{\Bbb
Q}_2} ``{\underset\tilde {}\to \nu_0},{\underset\tilde {}\to
\nu_1},{\underset\tilde {}\to \nu_2} \in \dsize \prod_{\ell < \omega}
\lambda_\ell$ form a counterexample that is: ${\underset\tilde {}\to
\nu_0} E {\underset\tilde {}\to \nu_1} \wedge 
{\underset\tilde {}\to \nu_1} E {\underset\tilde {}\to \nu_2} 
\wedge \neg \nu_0 E {\underset\tilde {}\to \nu_1}$ or $\neg
{\underset\tilde {}\to \nu_0} E {\underset\tilde {}\to \nu_0}$ or
${\underset\tilde {}\to \nu_0} E {\underset\tilde {}\to \nu_1} \wedge
\neg {\underset\tilde {}\to \nu_1} E {\underset\tilde {}\to \nu_0}$ or
${\underset\tilde {}\to \nu_0} E {\underset\tilde {}\to \nu_1}
\wedge (\exists !n)({\underset\tilde {}\to \nu_0}(n) \ne 
{\underset\tilde {}\to \nu_1}(n))"$. \nl
Choose $\chi$ large enough and $\bar N = \langle N_n:n < \omega
\rangle,N$ such that:
\mr
\item "{$\circledast^\chi_{\bar N}(i)$}"  $N_n
\prec_{L_{\lambda^+_n,\lambda^+_n}} ({\Cal H}(\chi),\in)$ and $\|N_n\|
= 2^{\lambda_n}$ and $\{p^*,E,{\underset\tilde {}\to
\nu_0},{\underset\tilde {}\to \nu_1},
{\underset\tilde {}\to \nu_2},N_0,\dotsc,N_{n-1}\}$ belong to $N_n$
\sn
\item "{$(ii)$}"  $N = \dbcu_{n < \omega} N_n$ so $N \prec ({\Cal
H}(\chi),\in)$.
\ermn
Now we choose $p_n$ by induction on $n < \omega$ such that:
\mr
\item "{$(*)(i)$}"  $p_0=p$,
\sn
\item "{$(ii)$}"  $p_n \in N_n,n(p_n) \ge n$
\sn
\item "{$(iii)$}"  $p_n \le p_{n+1}$
\sn
\item "{$(iv)$}"  if $\underset\tilde {}\to \tau \in N_n$ is a $\Bbb
Q_2$-name of an ordinal then for some $k_n(\underset\tilde {}\to \tau)
> n+1$ we have: if $p_{n+1} \le q$ and $n(q) \ge k_n({\underset\tilde {}\to
\tau_n})$ then $q$ forces a value to $\underset\tilde {}\to \tau$.
\ermn
This is possible by \scite{3.9}(2),(3).  Now let $G = \{q:q \in N \cap
\Bbb Q_2$ and $q \le p_n$ (or just $p_n \Vdash ``q \in \underset\tilde
{}\to G")$ for some $n\}$, it is a subset of $\Bbb Q^N_2$ generic over $N$.
(Why?  If $N \models ``{\Cal I} \subseteq \Bbb Q_2$ is dense" then ${\Cal I}
\subseteq \Bbb Q_2$ is dense and there is ${\Cal I}' \subseteq {\Cal
I}$, a maximal antichain of $\Bbb Q_2$ which belongs to $N$ hence to
some $N_n$; there is $g \in N_n$, a one to one function from ${\Cal
I}'$ onto $|{\Cal I}'|$, so it define a $\Bbb Q_2$-name
$\underset\tilde {}\to \tau,
\underset\tilde {}\to \tau(G) = \gamma \Leftrightarrow (\forall q)(q
\in {\Cal I}' \cap G \rightarrow f(q) = \gamma) \Leftrightarrow
(\exists q)(q \in {\Cal I}' \cap G \and f(q) = \gamma)$, 
so $k_n(\underset\tilde {}\to \tau) < \omega$ is well 
defined and so $p_{k_n(\underset\tilde {}\to \tau)}$
forces a value to $\underset\tilde {}\to \tau$ hence forces $q \in
\underset\tilde {}\to G$ for some $q \in {\Cal I}' \subseteq {\Cal
I}$, hence $q \in G$ so $G \cap {\Cal I} \ne \emptyset$ as required).
Now by straightforward absoluteness argument, ${\underset\tilde {}\to
\nu_0}[G],
{\underset\tilde {}\to \nu_1}[G],{\underset\tilde {}\to \nu_2}[G] \in
\dsize \prod_{\ell < \omega} \lambda_\ell$ give contradiction to an
assumption.  \hfill$\square_{\scite{3.11}}$
\enddemo
\bn
In fact
\demo{\stag{3.11A} Observation}  Assume
\mr
\item "{$(a)$}"  $\lambda^*$ is strong limit of cofinality
$\aleph_0,\lambda^* = \dsize \sum_{n < \omega} \lambda^*_n,\lambda^*_n <
\lambda^*_{n+1}$
\sn
\item "{$(b)$}"
{\roster
\itemitem{ $(i)$ }   $\Bbb Q$ is a forcing notion
\sn
\itemitem{ $(ii)$ }   $\le_{pr}$ is included in $\le_{\Bbb Q}$ and
\sn
\itemitem{ $(iii)$ }  $\bold n:\Bbb Q \rightarrow \omega$ is such
that, for each $n$ the set ${\Cal I}_n = \{p \in \Bbb Q:\bold n(p) \ge
n\}$ is a dense subset of $\Bbb Q$
\sn
\itemitem{ $(iv)$ }  for $p \in \Bbb Q,\{q \in \Bbb Q:p \le_{pr} q\}$ is
$\lambda^*_{n(p)}$-complete and
\sn
\itemitem{ $(v)$ }   $\Bbb Q$ has pure decidability for $\Bbb
Q$-names of truth values
\sn
\itemitem{ $(vi)$ }  if $p \in \Bbb Q$ and $\underset\tilde {}\to
\tau$ is a $\Bbb Q$-name of an ordinal, \ub{then} there are $m <
\omega$ and $q$ satisfying: $p \le_{pr} q$ and $q \le r \wedge m \le
n(r) \Rightarrow$ ($r$ forces a value to $\underset\tilde {}\to \tau$)
\endroster}
\item "{$(c)$}"  $N,\langle N_n:n < \omega \rangle$ as above, 
$\{\Bbb Q,\le,\le_{pr}\} \in N_0$.
\ermn
\ub{Then} there is $G \subseteq \Bbb Q^N$ generic over $N$ hence
${\Cal H}(\lambda)^{N[G]} = {\Cal H}(\lambda)$.
\enddemo
\bigskip

\demo{Proof}  Should be clear.
\enddemo
\bigskip

\definition{\stag{3.12} Definition/Claim}  Assume that $F$ is a permutation of
$(\dsize \prod_{\ell < n(*)} \lambda_\ell) \times (\dsize \prod_{\ell
< n(*)} \lambda_\ell)$. \nl
1) Let $\Bbb Q^{\ge n(*)}_2 = \{p \in \Bbb P_2:n(p) \ge n(*)\}$ and let
$\hat F$ be the following function from $\Bbb Q^{\ge n(*)}_2$ to $\Bbb
Q^{\ge n(*)}_2$

$$
\align
\hat F(p) = q \text{ iff } &n(q) = n(p) \\
  &(\eta^q_0 \restriction n(*),\eta^q_1 \restriction n(*)) =
F((\eta^p_0 \restriction n(*),\eta^p_1 \restriction n(*))) \\
  &\eta^q_0 \restriction [n(*),n(p)) = \eta^p_0 \restriction
[n(*),n(p)) \\
  &\eta^q_1 \restriction [n(*),n(p)) = \eta^q_1 \restriction
[n(*),n(p)) \\
  &\bar A^q = \bar A^p.
\endalign
$$
\mn
2) For $p \in \Bbb Q^{\ge n(*)}_2,\hat F(p)$ is well defined $\in \Bbb
Q^{\ge n(*)}_2$. \nl
3) $\hat F$ is a permutation of $\Bbb Q^{\ge n(*)}_2$ preserving
$\le,\le_{pr},\le_{pr,n},\le_{apr}$ (and their negations), and $F \mapsto
\hat F$ is a group homomorphism. \nl
4) If $G \subseteq \Bbb Q_2$ is generic over $\bold V$ then
$\hat{\bold F}(G) = \{r \in \Bbb Q_2$: for some $q \in G \cap \Bbb
Q^{\ge n(*)}_2$ we have $r \le \hat F(q)\}$ is a generic over $\bold V$ and
$\bold V[\hat F(G)] = \bold V[G]$ and even $N[\hat F(G)] = N[G]$ if
say $N \prec ({\Cal H}(\chi),\in),\Bbb Q_2 \in N,F \in N,\lambda
\subseteq N$.
\enddefinition
\bigskip

\demo{Proof}  Easy.
\enddemo
\bigskip
 
\proclaim{\stag{3.13} Claim}

$\Vdash_{{\Bbb Q}_2} ``\neg {\underset\tilde {}\to \eta_0} E
{\underset\tilde {}\to \eta_1}"$.
\endproclaim
\bigskip

\demo{Proof}  If not, let $p \in \Bbb Q_2$ be such that $p \Vdash ``{\underset\tilde
{}\to \eta_0} E {\underset\tilde {}\to \eta_1}"$.
Now by clause (f) of Definition \scite{3.2}(1), we can find $p_1$ such
that:
\mr
\item "{$(i)$}"  $\Bbb Q_2 \models p \le_{pr} p_1$
\sn
\item "{$(ii)$}"  if $n(p) \le n < \omega$ and $\beta \in s_1(A^{p_1}_n)$
then for some $\alpha,\gamma$ we have $(\alpha,\beta),(\beta,\gamma)
\in s_2(A^p_n)$.
\ermn
Let $G_1 \subseteq \Bbb Q_2$ be generic over $\bold V$ such that $p_1 \in G_1$ and
let $\eta_\ell = {\underset\tilde {}\to \eta_\ell}[G_1]$ so $\bold V[G_1]
\models \eta_0 E \eta_1$.  By \scite{3.11} in $\bold V[G_1],E$ is still an
equivalence relation satisfying clause (d) of \scite{3.1} and
$\eta_1(n) \in s_1(A^{p_1}_n)$.  We can
find $\alpha < \lambda_{n(p)}$ such that $\alpha^* < \eta_1(n^*)$ and
$(\alpha^*,\eta_1(n^*)) \in s_2(A^p_{n^*})$.  Let us define $\eta'_0
\in \dsize \prod_{n < \omega} \lambda_n$ by $\eta'_0(n)$ is $\alpha^*$
if $n=n^*$ and $\eta_0(n)$ otherwise. \nl
Now the pairs $(\eta_0 \restriction (n(*)+1),\eta_1 \restriction
(n(*)+1))$ and $(\eta'_0 \restriction (n(*))+1),\eta_1 \restriction
(n(*)+1))$ are from $(\dsize \prod_{n \le n(*)} \lambda_n) \times
(\dsize \prod_{n \le n(*)} \lambda_n)$, so there is a permutation $F$
of this set interchanging those two pairs and is the identity
otherwise.  Let $G_2 = \hat F(G_1)$.  Now by \scite{3.12}
\mr
\item "{$(*)_1$}"  $G_2$ is a generic subset of $\Bbb Q_2$ over $\bold
V$
\sn
\item "{$(*)_2$}"  $\bold V[G_2] = \bold V[G_1]$
\sn
\item "{$(*)_3$}"  ${\underset\tilde {}\to \eta_0}[G_2] = \eta'_0,
{\underset\tilde {}\to \eta_1}[G_2] = \eta_1$.
\ermn
By \scite{3.11} we have
\mr 
\item "{$(*)_4$}"  $\bold V[G_\ell] \models \neg \eta_0 E \eta'_0$.
\ermn
As $p \le p_1 \in G_1$, by the choice of $p$ clearly
\mr
\item "{$(*)_5$}"  $\bold V[G_1] \models ``\eta_0 E \eta_1"$.
\ermn
By the choice of $p_1$ and $(\alpha,\eta_1(n^*))$ clearly $p \le \hat
F(p_1) \in G_2$ so
\mr
\item "{$(*)_6$}"  $\bold V[G_2] \models ``{\underset\tilde {}\to \eta'_0}[G_2]E
{\underset\tilde {}\to \eta_1}[G_2]"$ hence $\bold V[G_\ell] \models ``\eta'_0 E
\eta_1"$.
\ermn
Now $(*)_4 + (*)_5 + (*)_6$ contradict \scite{3.11}.
\hfill$\square_{\scite{3.13}}$
\enddemo
\bigskip

\proclaim{\stag{3.14} Claim}  1) Fix $\chi > \lambda$ large enough and
choose $N_n \prec_{L_{\lambda_n,\lambda_n}} ({\Cal H}(\chi),\in)$ such
that $\|N_n\| = 2^{\lambda_n},\{E,p_n:n < \omega\},\{N_\ell:\ell <
n\}$ belongs to $N_n$, and let $N = \dbcu_{n < \omega} N_n$; (certainly
can be done).  \ub{Then} we can find $\langle \rho_\nu:\nu \in \dsize \prod_{\ell <
n} \mu_\ell,n < \omega \rangle$ and $N,\langle N_n:n < \omega
\rangle$
\mr
\item "{$(\alpha)$}"  $\rho_\nu \in \dsize \prod_{\ell < \ell g(\nu)} \lambda_\ell$ 
\sn
\item "{$(\beta)$}"  $\nu_1 \triangleleft \nu_2 \Rightarrow
\rho_{\nu_1} \triangleleft \rho_{\nu_2}$
\sn
\item "{$(\gamma)$}"  if $\nu_1,\nu_2 \in \dsize \prod_{\ell <n}
\lambda_\ell$ and $m \le k < n$ and $\nu_2(m) < \nu_2(m)$ then
$\eta_{\nu_1}(k) < \eta_{\nu_2}(m)$
\sn
\item "{$(\delta)$}"  if $\nu \in \dsize \prod_{i < \omega}
\lambda_\ell$ then $\rho_\nu =: \dbcu_{n < \omega} \rho_{\nu
\restriction n}$ is generic for $(N,\Bbb Q_1)$ 
\sn
\item "{$(\varepsilon)$}"  if $\nu_0,\nu_1 \in \dsize \prod_{\ell < \omega}
\lambda_\ell,\nu_0 <_{lex} \nu_1$ then $(\rho_{\nu_0},\rho_{\nu_1})$
is generic for $(N,\Bbb Q_2)$ hence
\sn
\item "{$(\zeta)$}" if $\nu_0 \ne \nu_1 \in \dsize \prod_{\ell <
\omega} \lambda_\ell$ then $\neg \eta_{\nu_0} E \eta_{\nu_1}$.
\ermn
2) Moreover, for some $p \in \Bbb Q_2,n(p) = 0$ and non-principal
ultrafitler $D$ on $\omega$ we have
\mr
\item "{$(*)$}"  if $\eta,\nu \in \dsize \prod_{n < \omega}
s_1(A^p_n)$ and $\eta/D \ne \nu/D$ then $\neg \eta E \nu$.
\endroster
\endproclaim
\bigskip

\demo{Proof}   Let $M_0 \prec_{L_{\aleph_1,\aleph_1}} N_0,\|M_0\| =
2^{\aleph_0}$, let
$N = \dbcu_{n < \omega} N_n$. \nl
As above we choose $p_n$ by induction on $n$ such that:
\mr
\widestnumber\item{$(iii)$}
\item "{$(i)$}"  $p_n \in \Bbb Q_n$
\sn
\item "{$(ii)$}"  $p_n \in N_n$
\sn
\item "{$(iii)$}"  $n(p_0) = 0$
\sn
\item "{$(iv)$}"  $p_n \le_{pr} p_{n+1}$
\sn
\item "{$(v)$}"  for every $\Bbb Q_2$-name of an ordinal
$\underset\tilde {}\to \tau \in N_n$, for some $k_n(\underset\tilde
{}\to \tau) \in [n,\omega)$ we have: if $\Bbb Q_2 \models ``p_{n+1} \le q"$
and $n(q) \ge k_n(\underset\tilde {}\to \tau)$ then $q$ forces a value to
$\underset\tilde {}\to \tau$
\sn
\item "{$(vi)$}"  if $\underset\tilde {}\to \tau \in M_0$ is a $\Bbb
Q_2$-name of a natural number then $p_0$ forces a value to it.
\ermn
We can find $p_\omega \in \Bbb Q_2$ such that $n < \omega \Rightarrow
p_n \le_{pr} p_\omega$ and we can find $p^*$ such that $p_\omega \le
p^*$ and $(\forall n)(\forall \beta)(\exists \alpha,\gamma)[\beta \in
s_1(A^{p^*}_n) \rightarrow (\alpha,\beta),(\beta,\gamma) \in
s_2(A^{p_\omega}_n)]$ and we shall show that $p^*$ is as required.
Now clearly
\mr
\item "{$\boxtimes$}"  if $\eta_0,\eta_1 \in \dsize \prod_{n < \omega}
\lambda_n$ and $(\forall \ell < 2)(\forall n < \omega)(\eta_\ell(n) \in
s_1(A^{p^*}_n)$ and for every $n < \omega$ large enough
$(\eta_0(n),\eta_1(n)) \in s_2(A^{p^*}_n)$ \ub{then}
{\roster
\itemitem{ $(a)$ }  for some subset $G$ of $\Bbb Q^N_2$ generic over $N$
we have ${\underset\tilde {}\to \eta_0}[G] = \eta_0,{\underset\tilde
{}\to \eta_1}[G] = \eta_1$
\sn
\itemitem{ $(b)$ }  $\neg \eta_0 E \eta_1$ \nl
[why? by \scite{3.12}.]
\endroster}
\ermn
This suffices for part (1): by clause (b) of Definition \scite{3.2}(1),
we can find $Y_n \subseteq \lambda_n$ of
cardinality $\mu^+_n$ (really otp$(Y_n) = \mu_n \times
\mu_{n-1} \times \ldots \times \mu_0$ is enough) such that for
any $\alpha < \beta$ from $Y,(\alpha,\beta) \in s_2(A^p_n)$.  Now we
can choose by induction on $n,\langle \rho_\nu:\nu \in \dsize
\prod_{\ell < n} \mu_\ell \rangle$ as required in
$(\alpha),(\beta),(\gamma)$ of \scite{3.14}(1), they are as required.
\sn
For $B \subseteq \omega$ let ${\underset\tilde {}\to \eta_B}$ be the
following $\Bbb Q_2$-name:

$$
\eta_B(n) \text{ is } {\underset\tilde {}\to \eta_0}(n) \text{ if } n
\in B \text{ and is } {\underset\tilde {}\to \eta_1}(n) \text{ if } n
\in \omega \backslash B.
$$
\mn
Clearly ${\underset\tilde {}\to \eta_B}$ is a $\Bbb Q_2$-name of a
member of ${}^\omega \lambda$ and ${\underset\tilde {}\to \eta_B} \in
M_0$ hence for $B_1,B_2 \subseteq \omega$ the following $\Bbb Q_2$-name of a
truth value, Truth Value$({\underset\tilde {}\to \eta_{B_1}} E
{\underset\tilde {}\to \eta_{B_2}})$, is decided by $p_0$, say it is
$\bold t(B_1,B_2)$.
\sn
Define a two place relation $E'$ on ${\Cal P}(\omega):B_2 E' B_2$ iff
$\bold t(B_1,B_2) =$ truth. \nl
Let $J = \{B \subseteq \omega:\bold t(\emptyset,B) =$ truth$\}$. \nl
Clearly
\mr
\item "{$(*)_0$}"  $E'$ is an equivalence relation on ${\Cal
P}(\omega)$
\sn
\item "{$(*)_1$}"  $\omega \notin J$, even $[n,\omega) \notin J$.
\ermn
Let $\alpha^1_n < \alpha^2_n < \alpha^3_n < \alpha^4_n$ be from $Y_n$
for $n < \omega$ and for $h \in {}^\omega \{1,2,3,4\}$ let $\nu_h \in \dsize
\prod_{n < \omega} \lambda_n$ be $\nu_h(n) = \alpha^{h(n)}_n$.  Easily
\mr
\item "{$(*)_2$}"  if $B_1,B_2,B_3,B_4 \subseteq \omega$ and $B_1
\backslash B_2 = B_3 \backslash B_4,B_2 \backslash B_1 = B_4
\backslash B_3$ then $\bold t(B_1,B_2) = \bold t(B_3,B_4)$ that is $B_1 E' B_2 \equiv
B_3 E B'_4$ 
\sn
\item "{$(*)_3$}"  if $B_1,B_2,B_3,B_4 \subseteq \omega$ and $B_2
\triangle B_2 = B_3 \triangle B_4$ then $B_1 E' B_2 \equiv B_3
E' B_1$ \nl
[why? this is the meaning of $(*)_2$]
\sn
\item "{$(*)_4$}"  for $B_1,B_2 \subseteq \omega,B_1 E' B_2$ iff $B_1
\triangle B_2 \in J$ \nl
[why?  use $(*)_2$ with $B_1,B_2,B_1 \triangle B_2,\emptyset$ here
standing for $B_1,B_2,B_3$ by there] 
\sn
\item "{$(*)_5$}"  $B_1 \subseteq B_2 \in J \Rightarrow B_2 \in J$ \nl
[why?  we can find $A_1,A-2,A_3 \subseteq \omega$ such that $(A_1
\triangle A_2) = B_2 = (A_2 \triangle A_3)$ and $B_1 = A_1 \triangle
A_3$ by the if part of $(*)_4$ we have $\bold t(A_1,A_2) =$ truth,
$\bold t(A_2,A_3) =$ truth, by the definition of $\bold t$ we have
$\bold t(A_1,A_3) =$ truth, hence by $(*)_4$ the only if part we have
$B_1 = A_1 \triangle A_3 \in J$ as required] 
\sn
\item "{$(*)_6$}"  $J$ is an ideal \nl
[why?  by $(*)_5$ it is enough to show that for disjoint $B_1,B_2 \in
J$ also $B = B_1 \cup B_2 \in J$.  Now we have $B \triangle B_1 = B_2$
hence by $(*)_4$ we have $B E' B_1$ and $B_1 \triangle \emptyset = B_1
\in J$ hence $B_1 E' \emptyset$ hence $B E' \emptyset$; but by $(*)_4$
this implies $B \in J$ as required.] 
\ermn
So by $(*)_6$ there is an ultrafilter $D$ on $\omega$ disjoint to $J$,
and by $(*)_1$ it is non-principal, so we have proved also part (2).
\hfill$\square_{\scite{3.14}} \quad \square_{\scite{3.1}}$
\enddemo
\newpage

\head {\S4 The countable cofinality case: negative results} \endhead  \resetall \sectno=4
\bigskip

\proclaim{\stag{d.1} Claim}  Assume
\mr
\item "{$(a)$}"  $\lambda > \text{ cf}(\lambda) = \aleph_0$ 
\sn
\item "{$(b)$}"   $(\forall \alpha < \lambda)[|\alpha|^{\aleph_0} < \lambda]$
\sn
\item "{$(c)$}"   there is an algebra ${\frak B}$ with universe
$\lambda$, no infinite free subset and with $< \lambda$ functions.
\ermn
\ub{Then} there is $E$ such that
\mr
\item "{$(\alpha)$}" $E$ is an equivalence relation on ${}^\omega \lambda$
\sn
\item "{$(\beta)$}"  $E$ is very nice (see Definition \scite{0.2})
\sn
\item "{$(\gamma)$}"  if $\eta,\nu \in {}^\omega \lambda$ and $\eta
=^* \nu$ (i.e. $(\exists^{< \aleph_0} n)(\eta(n) \ne \nu(n))$ then
$\eta E \nu \Leftrightarrow \eta = \nu$
\sn
\item "{$(\delta)$}"  $E$ has $\lambda$ equivalence classes.
\endroster
\endproclaim 
\bigskip

\remark{\stag{d.2} Remark}  1) We can replace ${}^\omega \lambda$ by the set of
increasing $\omega$-sequences or by $\dsize \prod_{n < \omega}
\lambda_n$ when $\lambda_n < \lambda_{n+1} < \lambda = \dsize \sum_{m
< \omega} \lambda_m$ or by $\{A \subseteq \lambda:(\forall
n)(\exists!\alpha)(\alpha \in A \and \dsize \sum_{\ell < n}
\lambda_\ell \le \alpha < \lambda_n)$. \nl
2)  We can omit clause (b) if we weaken clause $(\gamma)$.  We can 
imitate \scite{3.3}, \scite{3.4}, \scite{3.5}.
\endremark 
\bigskip

\demo{Proof}  Without loss of generality ${\frak B}$ has $\aleph_0$
function and the individual constants $\{\alpha:\alpha < \lambda_0\}$
where $\lambda_0 < \lambda$.  

We define $E_0$ on ${}^\omega \lambda$ by

$$
\align
\eta E_0 \nu \text{ \ub{iff}}: &\text{ if } n <
\omega,\sigma(x_1,\dotsc,x_{n-1}) \text{ a } {\frak B} \text{-term,
and} \\
  &k,k_1,\dotsc,k_n < \omega \text{ then} \\
  &\eta(k) = \sigma(\eta(k_1),\eta(k_2),\dotsc,\eta(k_n)) \Leftrightarrow
\nu(k) = \sigma(\nu(k_1),\nu(k_2),\dotsc,\nu(k_n)).
\endalign
$$
\mn
So $E_0$ is an equivalence relation with $\le \lambda^{\aleph_0}_0 <
\lambda$ equivalence classes.  For $\eta \in {}^\omega \lambda$ let
$A(\eta) = \{k:\text{ for some } k^* < \omega$ there are no $n <
\omega,k_1,\dotsc,k_n \in [k^*,\omega)$ and ${\frak B}$-term
$\sigma(x_1,\dotsc,x_n)$ such that $\eta(k) =
\sigma(\eta(k_1),\dotsc,\eta(k_n))\}$. \nl
Lastly, we define $E_1$ by

$$
\eta E_1 \nu \text{ iff } \eta E_0 \nu \and \eta \restriction A(\eta)
= \nu \restriction A(\nu).  
$$ 
\mn
The rest is as in \S3.  \hfill$\square_{\scite{d.1}}$
\enddemo
\bigskip

\proclaim{\stag{d.3} Claim}  1) In \scite{d.1} we can demand
\mr
\item "{$(\delta)$}"  for each $\eta \in {}^\omega
\lambda,\eta/J^{bd}_\omega$ is a set of representatives of $E$.
\ermn
2) We can weaken in \scite{d.1} assumption (b) to
\mr
\item "{$(b)'$}"  $(\aleph_0 + |\tau({\frak B})|)^{\aleph_0} < \lambda$.
\ermn
3) If in \scite {d.1} we change clause $(\gamma)$ in the conclusion to
$(\gamma)^-$ below, we can omit clause (b) of the assumption
\mr
\item "{$(\gamma)^-$}"  for every equivalence class $X$ of $E$ for
every $\eta \in {}^\omega \lambda$ and for some
$\alpha,\eta_{n,\alpha} \in X$ where $\eta_{\alpha,n} \in {}^\omega
\lambda$ is: $\eta_{\alpha,n}(\ell) = \alpha$ if $\ell = n$ and
$\eta_{\alpha,n}(\ell) = \eta(\ell)$ otherwise.
\endroster
\endproclaim
\bigskip

\demo{Proof}  1) We imitate \scite{3.3} only letting $\Xi_\eta =
\bigl\{ \{ \langle k,k_1,\dotsc,k_n,\sigma \rangle:\nu(k) =
\sigma(k_1,\dotsc,k_n)\}:\nu \in {}^\omega \lambda,\nu/J^{bd}_\omega =
\eta/J^{bd}_\omega \bigr\}$. \nl
2) The same proof. \nl
3) For $\eta \in {}^\omega \lambda$ let $\bold n(\eta) < \omega$ be
the minimal $n \in [\bold n(\eta),\omega) \Rightarrow c \ell_{\frak
B}\{\eta(\ell):\ell \in [n,\omega)\} = c \ell_{\frak
B}\{\eta(\ell):\ell \in [\bold n(\eta),\omega)\}$.  Let $K$ be an
additive group with universe $\lambda,K_1$ a subgroup, $|K_1| =
\lambda,(K:K_1) = \lambda$ and $\eta E \nu$ iff $\dsize \sum_{n <
\bold n(\eta)} \eta(n) = \dsize \sum_{n < n(\nu)} \nu(\eta)$ mod
$K_1$. \nl
${{}}$  \hfill$\square_{\scite{d.3}}$
\enddemo
\bn
\ub{Question}:  What about having $\sigma \in (\lambda,2^\lambda)$
equivalence classes?
\newpage

\head {\S5 On $r_p(\text{Ext}(G,\Bbb Z)$} \endhead  \resetall \sectno=5
\bigskip

\definition{\stag{e.1} Definition}  For an abelian group $G$ and prime
$p$ let $r_p(G)$ be the rank of $G/pG$ as a vector space over $\Bbb
Z/p \Bbb Z$.

Instead using a definition of the abelian group Ext$(G,\Bbb Z)$, for
$G$ torsion free abelian group we quote (see \cite{Fu}).
\enddefinition
\bigskip

\proclaim{\stag{e.2} Claim}  For a torsion free abelian group $G$ and
prime $p,r_p(\text{\rm Ext\/}(g,\Bbb Z))$ is the rank of Hom$(G,\Bbb
Z/pZ)/(\text{Hom}(G,\Bbb Z)/p \Bbb Z)$ where Hom$(G,\Bbb Z/p \Bbb Z)$ is
the abelian group of homomorphisms from $G$ to $\Bbb Z/p \Bbb Z$,
Hom$(\Bbb Z,G)/p \Bbb Z$ is the abelian group of homomorphism $h$ from $G$
to $\Bbb Z/p \Bbb Z$ such that for some homomorphism $g$ from $G$ to
$\Bbb Z$ we have $x \in G \Rightarrow g(x)/\Bbb Z \in h(x)$.
\endproclaim
\bn
More generally (see \cite{Sh:664})
\definition{\stag{e.3} Definition}  1) We say ${\Cal Y} = (\bar A,\bar
K,\bar G,\bar D,\bar g^*)$ is a $\lambda$-system if
\mr
\item "{$(A)$}"  $\bar A = \langle A_i:i \le \lambda \rangle$ is an
increasing sequence of sets, $A = A_\lambda = \cup \{A_i:i < \lambda\}$
\sn
\item "{$(B)$}"  $\bar K = \langle K_t:t \in A \rangle$ is a sequence
of finite groups
\sn
\item "{$(C)$}"  $\bar G = \langle G_i:i \le \lambda \rangle$ is a
sequence of groups, $G_i \subseteq \dsize \prod_{t \in A_i} K_t$, each
$G_i$ is closed and $i < j \le \lambda \Rightarrow G_i = \{g
\restriction A_i:g \in G_j\}$ and \nl
$G_\lambda = \{g \in \dsize \prod_{t \in A_\lambda} K_t:(\forall i <
\lambda)(g \restriction A_i \in G_i)\}$ 
\sn
\item "{$(D)$}"  $\bar D = \langle D_\delta:\delta \le \lambda$ (a
limit ordinal)$\rangle,D_\delta$  an ultrafilter on $\delta$ such that
$\alpha < \delta \Rightarrow [\alpha,\delta) \in D_\delta$
\sn
\item "{$(E)$}"  $\bar g^* = \langle g^*_i:i < \lambda \rangle,g^*_i
\in G_\lambda$ and $g^*_i \restriction A_i = e_{G_i} = \langle
e_{K_t}:t \in A_i \rangle$.
\ermn
Of course, formally we should write $A^{\Cal Y}_i,K^{\Cal Y}_t,G^{\Cal
Y}_\delta,D^{\Cal Y}_\delta,g^{\Cal Y}_i$, eta., if clear from the
context we shall not write this. \nl
2) Let ${\Cal Y}^-$ be the same omitting $D_\lambda$ and we call it a
lean $\lambda$-system.
\enddefinition
\bigskip

We can deduce Sageev Shelah \cite{SgSh:148} result (if $|G| = \lambda$
is weakly compact $(> \aleph_0)$ and $p$ is prime, \ub{then}
$r_p(\text{Ext}(G,\Bbb Z)) \ge \lambda \Rightarrow
r_p(\text{Ext}(G,\Bbb Z)) = 2^\lambda)$.  For this note
\proclaim{\stag{e.8} Claim}  1)  Assume
\mr
\item "{$(a)$}"  ${\Cal Y}$ is a $\lambda$-system
\sn
\item "{$(b)$}"  $\bar H = \langle H_i:i < \lambda \rangle$ is a
sequence of groups, $\bar \pi = \langle \pi_{i,j}:i < j < \lambda
\rangle$, \nl
$\pi_{i,j} \in \text{ Hom}(H_j,H_i)$, commuting
\sn
\item "{$(c)$}"  $\bar h = \langle h_i:i < \lambda \rangle,h_i \in
\text{ Hom}(H_i,G^{\Cal Y}_i)$, and $i < j < \lambda \and x \in
H_j \Rightarrow (h_j(x)) \restriction A_i = h_i(\pi_{i,j}(x))$
\sn
\item "{$(d)$}"  $H_\lambda,\pi_{i,\lambda} \, (i < \lambda)$ form
the inverse limit of $\langle H_i,\pi_{i,j}:i < j < \lambda \rangle$,
and $h = h_\lambda$ the inverse limit of $\langle h_i:i < \lambda
\rangle$
\sn
\item "{$(e)$}"  $E_h$ is the following 2-place relation on $G_\lambda:f_1 E f_2
\Leftrightarrow f_1 f^{-1}_2 \in \text{ Rang}(h)$.
\ermn
\ub{Then}
\mr
\item "{$(\alpha)$}"  $h \in \text{ Hom}(H_\lambda,G_\lambda)$
\sn
\item "{$(\beta)$}"   if $(\forall i < \lambda)(|A_i| \le \lambda \and
|H_i| \le \lambda)$, then $E_h$ is a $\Sigma^1_1[\lambda]$-equivalence
relation
\sn
\item "{$(\gamma)$}"  if $(\forall i < \lambda)(|A_i| < \lambda \and
|H_i| < \lambda)$ and $\lambda$ is weakly compact uncountable,
\ub{then} the 2-place relation $E = E_h$ on $G_\lambda$ defined by $f_1
E f_2 \Leftrightarrow f_1 f^{-1}_2 \in \text{ Rang}(h)$ is very nice
\sn
\item "{$(\delta)$}"  under $(\gamma)$'s assumptions, if
$(G:\text{Rang}(h)) \ge \lambda$ then $(G:\text{Rang}(h)) = 2^\lambda$.
\ermn
2) If for $\varepsilon < \varepsilon(*) \le \lambda$ we have $\langle
H^\varepsilon_i:i < \lambda \rangle,\langle \pi^\varepsilon_{i,j}:i <
j < \lambda \rangle,\langle h^\varepsilon_i:i \le \lambda \rangle$ as
above, and for every $\alpha < \lambda$ there are $f^\alpha_i \in
G_\lambda$ (for $i < \alpha$) such that $\neg(f^\alpha_i
E_{h^\varepsilon_\lambda} f^\alpha_j)$ for $i < j < \alpha$, \ub{then}
there are $f_i \in G$ for $i < 2^\lambda$ such that $i < j < 2^\lambda
\and \varepsilon < \varepsilon^* \Rightarrow \neg(f_i E_{h_\varepsilon} f_j)$.
\endproclaim 
\bigskip

\demo{Proof}  Straight.  (The main point is in clause $(\gamma)$ of
\scite{e.8}(1), the point is that if $f \in G_\lambda \backslash
\text{ Rang}(h_\lambda)$ then for some $i < \lambda$ we have $\pi_{i,\lambda}(f)
\in G_i \backslash \text{ Rang}(h_i)$ by the tree property of
$\lambda$).  \hfill$\square_{\scite{e.8}}$
\enddemo
\bigskip

\proclaim{\stag{g.12} Claim}  Assume
\mr
\item "{$(A)(a)$}"  $\lambda$ is a strong limit cardinal and $\theta$
is a compact cardinal $< \lambda$
\sn
\item "{$(b)$}"  $K_i$ is a group for $i < \lambda$
\sn
\item "{$(c)$}"  $I$ is a directed partial order, $t \in I \Rightarrow A(t)
\subseteq \lambda$ and $\dbcu_{t \in I} A(t) = \lambda$
\sn
\item "{$(d)$}"  for $t \in I,G_t$ is a subgroup of $\dsize \prod_{i
\in u} \{K_i:i \in A(t)\}$
\sn
\item "{$(e)$}"  for $s \le t$ from $I$ we have $A(s) \subseteq A(t)$ and
$f \in G_t \Rightarrow f \restriction A(s) \in G_s$
\sn
\item "{$(f)$}"  $G_\infty$ is the inverse limit
\sn
\item "{$(B)(a)$}"  $\varepsilon(*) \le \lambda$
\sn
\item "{$(b)$}"  for $\varepsilon < \varepsilon(*),\langle
H^\varepsilon_u,\pi^\varepsilon_{u,w}:u \le w$ from $I\rangle$
is an inversely directed system of groups
\sn
\item "{$(c)$}"  $h^\varepsilon_u \in \text{
Hom}(H^\varepsilon_u,G_u)$ for $u \in I,\varepsilon < \varepsilon(*)$
\sn
\item "{$(d)$}"
$H^\varepsilon_\infty,h^\varepsilon,h^\varepsilon_{\infty,u}$ are the
limit of the inverse system
\sn
\item "{$(e)$}"  $E_\varepsilon$ is the equivalence relation on
$G_\infty:f E_\varepsilon g 
\Leftrightarrow fg^{-1} \in \text{ Rang}(h^\varepsilon_\infty)$
\endroster
\sn
$(C) \qquad$  for every $\mu < \lambda$ we can find $\langle
f^\mu_\alpha:\alpha < \mu \rangle$ from $G_\infty$ such that \nl

\hskip30pt $\varepsilon < \mu \cap \varepsilon(*) \and \alpha < \beta \Rightarrow
\neg(f^\mu_\alpha E_\varepsilon f^\mu_\beta)$
\sn
$(D) \qquad \theta$ is $> \underset {i < \lambda}{}\to \sup|K_i| +
\underset {t \in I}{}\to \sup|A(t)|$ and also 
$\sup_{t \in I,\varepsilon < \varepsilon(*)} |H^\varepsilon_t|$.
\mn
\ub{Then} there are $f_\alpha \in G$ for $\alpha < 2^\lambda$ such
that $\varepsilon < \varepsilon(*) \and \alpha < \beta = 2^\lambda
\Rightarrow \neg(f_\alpha E_\varepsilon f_\beta)$.
\endproclaim
\bigskip

\demo{Proof}  Let $\kappa = \text{ cf}(\lambda),\langle \lambda_i:i <
\kappa \rangle$ be increasing with limit $\lambda$.  We can choose by
induction on $i < \kappa,I_i,A_i$ such that
\mr
\item "{$(\alpha)$}"  $A_i \subseteq \lambda,|A_i| < \lambda,j < i
\Rightarrow A_j \subseteq A_i,\lambda_i \subseteq A_i$
\sn
\item "{$(\beta)$}"  $I_i \subseteq I$ is directed, $|I_i| < \lambda,j
< i \Rightarrow I_j \subseteq I_i$ and $t \in I_i \Rightarrow A(t) \subseteq
A_i$
\sn
\item "{$(\gamma)$}"  if we restrict ourselves to $A_i,I_i$, there is
$\langle f^i_\alpha:\alpha < \mu_i \rangle$, such that $f^i_\alpha \in
G^{I_i}_\infty = \text{ Lim}_{I_i} \langle G_u,f_{u,w}:u \le w$ from
$I_i\rangle$ and $\varepsilon < \mu_i \cap \varepsilon(*) \and \alpha
< \mu_i \Rightarrow \neg(f^i_\alpha E^{I_i} f^i_\beta)$.
\ermn
Now we can apply \cite[\S3]{Sh:664}.  \hfill$\square_{\scite{g.12}}$
\enddemo
\bigskip

\proclaim{\stag{g.13} Claim}   Assume
\mr
\item "{$(a)$}"   $\lambda > \text{ cf}(\lambda) = \kappa$, and
$\kappa$ is a measurable cardinal, say $D$ a normal ultrafilter on $\kappa$
\sn
\item "{$(b)$}"  $G$ is a torsion free abelian group
\sn
\item "{$(c)$}"  $|G| = \lambda$
\sn
\item "{$(d)$}"  $p$ is a prime number.
\ermn
1) If $r_p(\text{Ext}(G,\Bbb Z)) \ge \lambda$ and $\lambda =
\lambda^{< \kappa}$ then $r_p(\text{Ext}(G,\Bbb Z))
\ge \lambda^\kappa$. \nl
2) Let $\lambda = \dsize \sum_{i < \kappa} \lambda_i$; then
$r_p(\text{Ext}(G,\Bbb Z))$ is $\dsize \prod_{i < \kappa} \mu_i/D$ for
some $\mu_i \le 2^{\lambda_i}$.
\endproclaim
\bigskip

\demo{Proof}  1) Let $\langle G_i:i < \kappa \rangle$ be an increasing
sequence of pure subgroups of $G$ with union $G$ satisfying $i < \kappa \Rightarrow
|G_i| < \lambda$.  Now
\mr
\item "{$(*)$}"  if $g \in \text{ Hom}(G,\Bbb Z/p \Bbb Z)$ and $i <
\kappa \Rightarrow g \restriction G_i \in \text{ Hom}(G,\Bbb Z)/p \Bbb
Z$ then $g \in \text{ Hom}(G,\Bbb Z)/p \Bbb Z$. \nl
Why?  Let $g \restriction G_i = h_i/p \Bbb Z$ where $h_i \in \text{
Hom}(G,\Bbb Z)$ and let $h$ a function from $G$ to $\Bbb Z$ be defined
as $h(x)=n \Leftrightarrow \{i < \kappa:h_i(x)=n\} \in D$.  Clearly $h
\in \text{ Hom}(G,\Bbb Z)$ and $g = h/p \Bbb Z$, as required.]
\ermn
The rest should be clear. \nl
2) As in part (1), letting $\mu_i = r_p(\text{Ext}(G_i,\Bbb Z))$.
\hfill$\square_{\scite{g.13}}$
\enddemo
\bigskip

\demo{\stag{g.14} Conclusion}  If $\lambda$ is a strong limit and above a compact
cardinal and $G$ is a torsion free abelian group and $p$ is a prime
then $r_p(\text{Ext}(G,\Bbb Z)) \ge \lambda \Rightarrow
r_p(\text{Ext}(G,\Bbb Z)) = 2^\lambda$.
\enddemo
\bigskip

\demo{Proof}  By \scite{g.13}.
\enddemo
\bigskip

\remark{\stag{g.15} Remark} So for $\lambda$ strong limit singular the problem of
the existence of $G$ such that $|G| = \lambda,r_p(\text{Ext}(G,\Bbb
Z)) = \lambda$ is not similar to the problem of the existence of $M$
such that $\|M\| = \lambda$, nu$(M) = \lambda$. 
\endremark


\shlhetal

\newpage
    
REFERENCES.  
\bibliographystyle{lit-plain}
\bibliography{lista,listb,listx,listf,liste}

\enddocument

\bye